\documentclass[a4paper,11pt]{article}
\pdfoutput=1 

\usepackage{graphicx}
\usepackage{tikz}
\usepackage[inline]{enumitem}
\usepackage{amssymb}
\usepackage{amsmath}
\usepackage{hyperref}
\usepackage{cleveref} 

\newcommand{\comments}[1]{#1} 

\usepackage{amsfonts}
\usepackage{rotating}
\usepackage{lscape}
\usepackage{verbatim}

\usepackage{setspace}
\usepackage{color}
\usepackage{caption}
\usepackage{url}
\usepackage{rotating}
\allowdisplaybreaks

\DeclareMathOperator*{\argmax}{arg\,max}
\DeclareMathOperator*{\argmin}{arg\,min}
\usepackage{csquotes}

\newcommand{\cg}[1]{#1}
\newcommand{\js}[1]{#1}

\usepackage[margin=0.5in]{geometry}
\newcommand{\CGrev}[1]{#1}
\newcommand{\JSrev}[1]{#1}

\newcommand{\revtwoCG}[1]{#1}
\newcommand{\revtwoJS}[1]{#1}

\newcommand{\revtwoCGblue}[1]{#1}

\begin{document}

\title{Optimization Problems for Machine Learning: A Survey}
\author{Claudio Gambella $^1$, Bissan Ghaddar $^2$, Joe Naoum-Sawaya$^2$}
\date{$^1$ IBM Research Ireland, Mulhuddart, Dublin 15, Ireland, $^2$ Ivey Business School, University of Western Ontario, London, Ontario N6G 0N1, Canada }

\maketitle

\begin{abstract}
\JSrev{This paper surveys the machine learning literature and presents in an optimization framework several commonly used machine learning approaches.  Particularly, mathematical optimization models are presented for regression, classification, clustering, deep learning, and adversarial learning, as well as new emerging applications in machine teaching, empirical model learning, and Bayesian network structure learning.} Such models can benefit from the advancement of numerical optimization techniques which have already played a distinctive role in several machine learning settings. The strengths and the shortcomings of these models are discussed and potential research directions and open problems are highlighted. 
\end{abstract}


\tableofcontents

\section{Introduction}
\label{Intro}
The pursuit to create intelligent machines that can match and potentially rival humans in reasoning and making intelligent decisions goes back to at least the early days of the development of digital computing in the late 1950s \cite{solomonoff1957inductive}. The goal is to enable machines to perform cognitive functions by learning from past experiences and then solving complex problems under conditions that are varying from past observations. Fueled by the exponential growth in computing power and data collection coupled with the widespread of practical applications, machine learning is nowadays a field of strategic importance.

\subsection{\JSrev{Machine Learning Basics}}\label{sec:basics}
Broadly speaking, machine learning relies on learning a model that returns the correct output given a certain input. \revtwoCG{The inputs, i.e.\revtwoCGblue{,} predictor measurements, are typically values that represent the parameters that define a problem, while the output, i.e.\revtwoCGblue{,} response, is a value that represents the solution.} 

Machine learning models fall into two categories: supervised and unsupervised learning \cite{friedman2001elements, james2013introduction}. In supervised learning, a response measurement is available for each observation of predictor measurements and the aim is to fit a model that accurately predicts the response of future observations. \JSrev{More specifically, in supervised learning, values of both the input $x$ and the corresponding output $y$ are available and the objective is to learn a function $f$ that approximates with a reasonable margin of error the relationship between the input and the corresponding output. The accuracy of a prediction is evaluated using a loss function $\mathcal{L}(f(x),y)$ which computes a distance measure between the predicted output and the actual output. In a general setting, the best predictive model $f^*$ is the one that minimizes the risk 
\begin{align*}
 \mathbb{E}_p[\mathcal{L}(f(x),y)]=\int \int p(x,y) \mathcal{L}(f(x),y)dxdy \revtwoCG{,}
\end{align*}
where $p(x,y)$ is the probability of observing data point $(x,y)$ \cite{vapnik2013nature}. In practice $p(x,y)$ is unknown, however the assumption is that an independent and identically distributed sample of data points $(x_1,y_1),\dots,(x_n,y_n)$ forming the training dataset is given. Thus instead of minimizing the risk, the best predictive model $f^*$ is the one that minimizes the empirical risk such that 
\begin{align*}
 f^*=\mbox{arg min }\frac{1}{n}\sum_{i=1}^n \mathcal{L} (f(x_i),y_i). 
\end{align*}

When learning a model, a key aspect to consider is model complexity. Learning a highly complex model may lead to overfitting\revtwoCGblue{,} which refers to having a model that fits the training data very well but generalizes poorly to other data. The minimizer of the empirical risk will often lead to overfitting, and hence \revtwoCGblue{has} a limited generalization property. Furthermore, in practice the data may contain noisy and incorrect values, i.e.\revtwoCGblue{,} outliers, which impacts the value of the empirical risk and subsequently the accuracy of the learned model. Attempting to find a model that perfectly fits every data point in the dataset is \revtwoCGblue{thus} not desired\revtwoCGblue{,} since the predictive power of the model will be diminished when points that are far from typical are fitted. Usually\revtwoCGblue{,} the choice of $f$ is restricted to a family of functions $F$ such that 
\begin{align}
 f^*=\mbox{arg }\min_{f\in F}\frac{1}{n}\sum_{i=1}^n \mathcal{L} (f(x_i),y_i). \label{emprisk}
\end{align}
\CGrev{
The degree of model complexity is generally dictated by the nature and size of the training data. While simpler models are advised for small training datasets that do not uniformly cover the possible data ranges, complex models need large data sets to avoid overfitting.
}
}

\revtwoCGblue{On the other hand, in unsupervised learning,} response variables are not available and the goal of learning is to understand the underlying characteristics of the observations. \JSrev{Unsupervised learning thus attempts to learn from the distribution of the data the distinguishing features and the associations in the data. As such\revtwoCGblue{,} the main use-case for unsupervised learning is exploratory data analysis\revtwoCG{,} where the purpose is to segment and cluster the samples in order to extract insights. While with \revtwoCG{supervised learning there is a clear} measure of accuracy by evaluating the prediction to the known response, in unsupervised it is difficult to evaluate the validity of the derived structure.}

The fundamental theory of machine learning models and consequently their success can be largely attributed to research at the interface of computer science, statistics, and operations research. The relation between machine learning and operations research can be viewed along three dimensions: (a) machine learning applied to management science problems, (b) machine learning to solve optimization problems, (c)  machine learning problems formulated as optimization problems.

\subsection{Machine Learning and Operations Research}
Leveraging data in business decision making is  nowadays mainstream as any business in today's economy is instrumented for data collection and analysis. While the aim of machine learning is to generate reliable predictions, management science problems deal with optimal decision making. Thus\revtwoCGblue{,} methodological developments that can leverage data predictions for optimal decision making is an area of research that is critical for future business value \cite{bertsimas2014predictive, kraus2020deep, MORTENSON2015583}. 

Another area of research at the interface of machine learning and operations research is using machine learning to solve hard optimization problems and particularly $\mathcal{NP}$-hard integer constrained optimization \cite{BonLZ18, khalil2016learning, khalil2017learning, LodZ17, VACLAVIK2018}. In that domain, machine learning models are introduced to complement existing approaches that exploit combinatorial optimization through structure detection, branching, and heuristics. 

Lastly, the training of machine learning models can be naturally posed as an optimization problem with typical objectives that include optimizing training error, measure of fit, and cross-entropy \cite{boct2011optimization, bottou2018optimization, Cur17, wright2017optimization}. In fact, the widespread adoption of machine learning is \revtwoCGblue{in part} attributed to the development of efficient solution approaches for these optimization problems\revtwoCGblue{,} which enabled the training of machine learning models. As we review in this paper, the development of these optimization models has largely been concentrated in areas of computer science, statistics, and operations research\revtwoCGblue{. However,} diverging publication outlets, standards, and terminology persist.

\subsection{Aim and Scope}\label{sec:aim}
 \CGrev{
	The aim of this paper is to present machine learning as optimization problems. For that, in addition to publications in classical operations research journals, this paper surveys machine learning and artificial intelligence conferences and journals, such as the conference on Association for the Advancement of Artificial Intelligence and the International Conference on Machine Learning. Furthermore, since machine learning research has rapidly accelerated with many important papers still in the review process, this paper also surveys a considerable number of relevant papers that are available on the arXiv repository. \CGrev{This paper also complements the recent surveys of \cite{bottou2018optimization, Cur17, wright2017optimization} which described methodological developments for solving machine learning optimization problems; \cite{bengio2018machine, LodZ17} which discussed how machine learning advanced the solution approaches of mathematical programming; \revtwoCG{\cite{corne2012synergies, OLAFSSON20081429}} which described the interactions between operations research and data mining; \cite{bennett2006interplay} which surveyed solution approaches to machine learning models cast as continuous optimization problems; and \cite{song2019review} which provided an overview on the various levels of interaction between optimization and machine learning}.}
 
 This paper presents optimization models for regression, classification, clustering, and deep learning (including adversarial attacks), as well as new emerging paradigms such as machine teaching and empirical model learning. Additionally, this paper highlights the strengths and the shortcomings of the  models from a mathematical optimization perspective and discusses potential novel research directions. \CGrev{This is to foster efforts in mathematical programming for machine learning. While important criteria for contributions in operations research are the convergence guarantees, deviation to optimality and speed increments with respect to benchmarks, machine learning applications have a partly different set of goals, such as scalability, reasonable execution time and memory requirement, robustness and numerical stability and, most importantly, generalization  \cite{bennett2006interplay}. It is therefore common for mathematical programming approaches to sacrifice optimality \revtwoJS{(local or global)} and convergence guarantees to obtain better generalization \revtwoCG{properties}, by adopting strategies such as early stopping~\cite{poggio2017theory}. }

Following this introductory section, regression models are discussed in Section~\ref{sec:regression} while classification and clustering models are presented in Sections~\ref{sec:classification} and \ref{sec:clustering}, respectively. \revtwoCG{Linear dimension reduction methods are reviewed in Section~\ref{sec:dim-red}. Deep learning models are presented in Section~\ref{sec:deeplearning}, while models for adversarial learning are discussed in Section~\ref{sec:adv_learning}. New emerging paradigms that include machine teaching and empirical model learning are presented in Section~\ref{sec:emerging}.} 
Finally, conclusions are drawn in Section~\ref{sec:conclusion}.

\section{Regression Models}\label{sec:regression}

\subsection{Linear Regression}
\revtwoCG{}

\revtwoCG{Since the early era of statistics, linear regression models have been widely adopted  in supervised learning for predicting a quantitative response}. 
The central assumption is that the \revtwoCG{relationship between} the dependent variables (\textit{feature measurements}, or \textit{predictors}, or \textit{input vector}) \revtwoCG{and} the independent \CGrev{variable} (real-valued \textit{output}\revtwoCGblue{, or \textit{response}}) is representable with a linear function (\textit{regression function}) with a reasonable accuracy.
\revtwoCG{Linear regression models preserve considerable interest, given their simplicity, their extensive range of applications, and the ease of interpretability. In particular, machine learning interpretability, in its simplest form, is the ability to explain in a humanly understandable way the role of the inputs in the outcome \cite{doshi2017towards}.
} 

\JSrev{Linear regression aims to find a linear function $f$ that expresses the relation between an input vector $x$ of dimension $p$ and a real-valued output $f(x)$ such as
\begin{align}
f(x) =  \beta_0 + x^\top \beta \revtwoCG{,}
\label{lin_regr_model} 
\end{align}
where $\beta_0\in \mathbb{R}$ is the intercept of the regression line and $\beta \in \mathbb{R}^p$ is the vector of coefficients corresponding to each of the input variables.} \JSrev{\JSrev{In order to estimate the regression parameters $\beta_0$ and $\beta$, one needs a training set $(X,y)$ where $X\in \mathbb{R}^{n\times p}$ denotes $n$ training inputs $x_1, \dots, x_n$ and $y$ denotes $n$ training outputs where each $x_i \in \mathbb{R}^{p}$ is associated with the real-valued output $y_i$.  The objective is to minimize the empirical risk \eqref{emprisk}\revtwoCG{, in order to quantify via $\beta_j$ the association between predictor $X_j$ and the response, for each $j=1, \dots, p$.}} 

The most commonly used loss function for regression is the \textit{least squared estimate}, where fitting a regression model reduces to minimizing the residual sum of squares (RSS) between the labels and the predicted outputs\revtwoCGblue{,} such as}
\begin{align}
RSS(\beta) = \sum_{i=1}^{n} (y_i - \beta_0 - \sum_{j=1}^{p} x_{ij} \beta_j)^2.
\label{lin_regr_RSS}
\end{align}
The least squares estimate is known to have the smallest variance among all linear unbiased estimates, and has a closed form solution. However, this choice is not always ideal \revtwoCGblue{for fitting}, since it can yield a model with low prediction accuracy, due to a large variance, and often leads to a large number of non-zero regression coefficients (i.e.\revtwoCGblue{,} low interpretability). \CGrev{Shrinkage methods discussed in Section \ref{sec:shrink} and Linear Dimension Reduction discussed in Section \ref{sec:dim-red}} are alternatives to the least squared estimate. \revtwoJS{Forward or backward elimination are also commonly used approaches to perform variable selection and to avoid overfitting \cite{friedman2001elements}.} 

The process of gathering input data is often affected by noise, which can impact the accuracy of statistical learning methods. \CGrev{A model that takes into account the noise in the features of linear regression problems is presented in \cite{BERTSIMAS2018931}, which also investigates the relationship between regularization and robustness to noise.} The noise is assumed to vary in an uncertainty set $\mathcal{U} \in \mathbb{R}^{n \times p}$, and the \CGrev{learner} adopts the robust \revtwoJS{perspective}: 
\begin{align}
\min_{\beta_0, \beta } \max_{\Delta \in \mathcal{U}} & \
g(y-\beta_0 - (X+\Delta) \beta) \revtwoCG{,}
\label{robust_regr} 
\end{align}
where $g$ is a convex function that measures the residuals (e.g., a norm function). The characterization of the uncertainty set $\mathcal{U}$ directly influences the complexity of problem \eqref{robust_regr}. 

The design of high-quality linear regression models requires several desirable properties, which are often conflicting and not simultaneously implementable. A fitting procedure based on \CGrev{Mixed Integer} Quadratic Programming (MIQP) is presented \CGrev{in \cite{bertsimasORforum}} and takes into account sparsity, joint inclusion of subset of features (called selective sparsity), robustness to noisy data, stability against outliers, modeler expertise, statistical significance, and low global multicollinearity. \CGrev{Mixed Integer Programming} (MIP) models for regression and classification are also investigated in \cite{bertsimas2007CR}. The regression problem is modeled as an assignment of data points to groups with the same regression coefficients.

\CGrev{In order to speed up the fitting procedure and improve the interpretability of the regression model, irrelevant variables can be excluded via feature selection strategies. For example, feature selection is desired in case some regression variables are highly correlated. \revtwoCGblue{Multicollinearity can be} detected by the condition number of the correlation matrix or the variance influence factor (VIF) \cite{chatterjee2015regression}. To \revtwoCGblue{achieve feature selection in this case}, \cite{tamura2017best} introduces a mixed integer semidefinite programming formulation to eliminate multicollinearity by bounding the condition number. The approach requires to solve a single optimization problem, in contrast with the cutting plane algorithm of \cite{bertsimasORforum}. Alternatively, \cite{tamura2019mixed} proposes a mixed integer quadratic optimization formulation with an upper bound on VIF, which is a better-grounded statistical indicator for multicollinearity with respect to the condition number. 
}

\subsection{Shrinkage methods}\label{sec:shrink}
{Shrinkage methods} \CGrev{(also called {regularization methods})} seek to diminish the value of the regression coefficients. The aim is to obtain a more interpretable model (with less relevant features), at the price of introducing some bias in the model determination. A well-known shrinkage method is {Ridge} regression, where a \JSrev{$2$-norm} penalization on the regression coefficients is added to the loss function such that
\begin{align}
\mathcal{L}_{ridge}(\beta_0, \beta) = \sum_{i=1}^{n} (y_i - \beta_0 - \sum_{j=1}^{p} x_{ij} \beta_j)^2 + \lambda \sum_{j=1}^p \beta_j^2 \revtwoCG{,}
\label{lin_regr_ridge}
\end{align}
where $\lambda$ controls the magnitude of shrinkage.
 \\ \ \\ \ \\
 
 \CGrev{
 Another technique for regularization in regression is the \textit{lasso regression}, which penalizes the $1$-norm of the regression coefficients, and seeks to minimize the quantity
 \begin{align}
\mathcal{L}_{lasso}(\beta_0, \beta) = \sum_{i=1}^{n} (y_i - \beta_0 - \sum_{j=1}^{p} x_{ij} \beta_j)^2 + \lambda \sum_{j=1}^p \lvert \beta_j \rvert.
 \label{lin_regr_lasso}
 \end{align}
  When  $\lambda$ is sufficiently large, the $1$-norm penalty forces some of the coefficient estimates to be exactly equal to zero, hence the models produced by the lasso are more interpretable than those obtained via Ridge regression.
}

 \CGrev{
Ridge and lasso regression belong to a class of techniques to achieve \textit{sparse regression}. As discussed in
\cite{bertsimas2016best, bertsimas2017sparseregr},  sparse regression can be formulated as the best subset selection problem \cite{miller2002subset}\JSrev{\begin{align}
	\min & \ \revtwoCG{\frac{1}{2} \lVert y - \beta_0- X \beta \rVert_2^2} \label{sparse_regr}\\
	\text{s.t }& \lVert \beta \rVert_0 \leq k, \label{sparse_regr:norm0}\\
			& \revtwoCG{\beta_0 \in \mathbb{R}, \beta \in \mathbb{R}^p,}\label{sparse_regr:var}
	\end{align}}
where $k$ is an upper bound on the number of predictors with a non-zero regression coefficient, i.e., the predictors to select, and \revtwoCG{$||\beta||_0$ is the number of non-zero entries of $\beta$, which is commonly referred to as the \enquote{$0$-norm} (though is not technically a norm as it does not satisfy the homogeneity property)} .
Problem \eqref{sparse_regr}--\eqref{sparse_regr:var} is $\mathcal{NP}$-hard\revtwoCGblue{, as proven in \cite{natarajan1995sparse}}. The recent work of \cite{bertsimas2016best} demonstrated that the best subset selection can be solved to near-optimal solutions using optimization techniques for values of $p$ in the hundreds or thousands. Specifically, by introducing the binary variables \JSrev{$s \in \{0,1\}^p$, the sparse regression problem can be transformed into the \revtwoCG{MIQP} formulation
	\begin{align}
		\min & \revtwoCGblue{\  \frac{1}{2} \lVert y - \beta_0 - X\beta \rVert_2^2} \label{sparse_regr_MIO} \\
	\text{s.t }& -M s_j \leq \beta_j \leq M s_j \quad \forall j = 1, \dots, p\revtwoCG{,} \label{sparse_regr:bigM}\\
	&  \sum_{j=1}^p s_j \leq k\revtwoCG{,}\\
		& \revtwoCG{\beta_0 \in \mathbb{R}, \beta \in \mathbb{R}^p,} \\
	& s \in \{0,1\}^p\revtwoCG{,} \label{sparse_regr_MIO:s}
	\end{align}
	where $M$ is a large constant, $M \geq \lVert \beta \rVert_{\infty}$.} Since the choice of the data
dependent constant $M$ largely affects the strength of the MIQP formulation, alternative formulations based on Specially Ordered Sets Type I can be devised \cite{claassen2007application}. 

As discussed in \cite{hastie2017extended}, the prediction accuracy of best subset selection is however highly dependent on the noise present in the input dataset, and it is not possible to establish a dominance relationship over lasso regression and forward stepwise selection \cite{efron2004least}.
In order to limit the effect of noise in the input data, \revtwoCGblue{make the model more robust,} and to avoid numerical issues, \cite{bertsimas2017sparseregr} introduces the Tikhonov regularization term
$\frac{1}{2\Lambda} \lVert \beta \rVert_2^2$
with weight $\Lambda > 0$  into the objective function of problem \eqref{sparse_regr_MIO}--\eqref{sparse_regr_MIO:s}, which is then solved using a cutting plane approach.}

 \CGrev{
The task of finding a linear model to express the relationship between regressors and \revtwoCG{output} is a particular case of selecting the hyperplane that minimizes a measure of the deviation of the data with respect to the induced linear form. 
As presented in \cite{blanco2018locating}, locating a hyperplane \revtwoCG{$\beta_0 + x^T \beta = 0, \beta_0 \in \mathbb{R} , \beta \in \mathbb{R}^{p}$ to fit a set of points $x_i \in \mathbb{R}^{p}, i = 1, \dots, n$, consists of finding $\hat{\beta_0}, \hat{\beta} \in \argmin_{\beta_0, \beta} \phi(\epsilon(\beta_0, \beta))$, where $\epsilon(\beta_0, \beta)=\epsilon_{\{x_1, \dots, x_n\}} (\beta_0, \beta)$} is a mapping to the residuals of the points on the hyperplane (according to a distance measure in $\mathbb{R}^p$), and $\phi$ is an aggregation function on the residuals (e.g., residual sum of squares, least absolute deviation \cite{edgeworth1887observations}). If the number of points $n$ is much smaller than the dimension $p$ of the space, feature selection strategies can be applied \cite{bertsimas2016best, miyashiro2015mixed}.
We note that hyperplane fitting is a variant of facility location problems \cite{diaz2004continuous, schobel1998locating}. 
}

\CGrev{\subsection{Regression Models Beyond Linearity}\label{sec:non-lin-regr}
\revtwoCG{A natural extension of linear regression models is to consider nonlinear terms}, which may capture complex relationships between regressors and predictors.  Nonlinear regression models include, among others,  polynomial regression, exponential regression, step functions, regression splines, smoothing splines and local regression \cite{friedman2001elements, james2013introduction}. Alternatively, the Generalized Additive Models (GAMs) \cite{GAMs} maintain the additivity of the original predictors $X_1, \dots, X_p$ and
the relationship between each feature and the response $y$ is expressed using nonlinear functions $f_j(X_j)$ such as
\begin{align}
y  =  \beta_0 + \sum_{j=1}^{p} f_j(X_j).
\label{GAM_model} 
\end{align}

GAMs may increase the flexibility and accuracy of the predictions with respect to linear models, while maintaining a certain level of interpretability of the predictors.  However, one limitation is given by the assumption of additivity of the features. To further increase the model flexibility, one could include predictors of the form $X_i \times X_j$, or consider non-parametric models, such as random forests and boosting.
It has been empirically observed that GAMs do not represent well problems where the number of observations is much larger than the number of predictors. In \cite{taylan2007new} the Generalized Additive Model Selection is introduced to fit sparse GAMs in high dimension with a penalized likelihood approach. The penalty term is derived from the fitting criterion for smoothing splines. Alternatively, \cite{chouldechova2015generalized} proposes to fit a constrained version of GAMs by solving a conic programming problem.}

\CGrev{As an intermediate model between linear and nonlinear relationships, compact and simple representations via piecewise affine models have been discussed in \cite{KESHVARI2018585}. Piecewise affine forms emerge as candidate models when the fitting function is known to be discontinuous \cite{ferrari2003clustering}, separable~\cite{DEFARIAS2008234}, or approximate to complex nonlinear expressions \cite{d2015mathematical, rovatti2014optimistic, vielma2010mixed}. Fitting piecewise affine models involves partitioning the domain $D$ of the input data into $K$ subdomains $D_i, i =1, \dots, K,$ and fitting for each subdomain an affine function $f_i: D_i \rightarrow \mathbb{R}$, in order to minimize a measure of the overall fitting error. To facilitate the fitting procedure, the domain is partitioned a priori (see $K$-hyperplane clustering in Section \ref{sec:k-hyp-clust}). Neglecting domain partitioning may lead to large fitting errors. In contrast, \cite{amaldi2016discrete} considers both aspects in determining piecewise affine models for piecewise linearly separable subdomains via a mixed integer linear programming formulation and a tailored heuristic. Mixed integer models are also proposed in \cite{TORIELLO201286}, however a partial knowledge of the subdomains is required. Alternatively, clustering techniques can be adopted for domain partitioning~\cite{ferrari2003clustering}.
}

\section{Classification}\label{sec:classification}
\JSrev{The task of classifying data is to decide the class membership of an \revtwoCG{unlabeled} data item $x$ based on the training dataset $(X,y)$ where each $x_i$ has a known class membership $y_i$.} A recent comparison of machine learning techniques for binary classification is found in \cite{BAUMANN20191041}. \revtwoCG{This section reviews the common binary and multiclass classification approaches} that include logistic regression, linear discriminant analysis, decision trees, and support vector machines.

\subsection{Logistic Regression}
 In  most  problem  domains,  there  is  no functional relationship $y=f(x)$ between $y$~and~$x$. In this case, the relationship between $x$ and $y$ has to be described more generally by a probability distribution $P(x,y)$ while assuming that the training contains independent samples from $P$. \revtwoJS{In this section, the label $y$ is assumed to be binary, i.e.\revtwoCGblue{,} $y \in \{0, 1\}$}. The optimal class membership decision is to choose the class label $y$ that maximizes the posterior distribution $P(y|x)$. Logistic regression calculates the class membership probability for one of the two categories in the dataset~as
  \begin{align*}
 &P(y=1|x,\beta_0,\beta)=h(x,\beta_0,\beta)= \frac{1}{1 + e^{\revtwoJS{-(\beta_0+\beta^{\top} x)}}},\\ 
 &P(y=0|x,\beta_0,\beta)=1-h(x,\beta_0,\beta).
\end{align*}
The decision boundary between the two binary classes is formed by a hyperplane whose equation is $\beta_0+\beta^\top x=0$. Points at this decision boundary have $P(1|x,\beta_0,\beta)=P(0|x,\beta_0,\beta)=0.5$. \revtwoCG{The parameters $\beta_0$ and $\beta$ are usually obtained by maximum-likelihood estimation \cite{dreiseitl2002logistic}
\begin{align*}
\max \ \Pi_{i=1}^n P(y_i|x_i,\beta_0,\beta) = \max \ \Pi_{i=1}^n(h(x_i,\beta_0,\beta))^{y_i}(1-h(x_i,\beta_0,\beta))^{1-y_i}, 
\end{align*}
which is equivalent to
\begin{align}
 \min-\sum_{i=1}^n (y_i\mbox{log} \ h(x_i,\beta_0,\beta) + (1-y_i)\mbox{log} (1-h(x_i,\beta_0,\beta))). \label{logit_training}
\end{align}
Problem \eqref{logit_training} is convex and differentiable and first order methods such as gradient descent as well as second order methods such as Newton's method can be applied to find a global optimal solution.}

To tune the logistic regression model and to avoid overfitting, variable selection can be performed where only the most relevant subsets of the $x$ variables are kept in the model \cite{friedman2001elements}. \JSrev{Heuristic approaches such as forward selection or backward elimination} can be applied to add or remove variables respectively, based on the statistical significance of each of the computed coefficients. Interaction terms can be also added to further complicate the model at the risk of overfitting the training data.


\subsection{Linear Discriminant Analysis}
Linear discriminant analysis (LDA) is an approach for classification and dimensionality reduction. It is often applied to data that contains a large number of features (such as image data) where reducing the number of features is necessary to obtain robust classification.  While LDA and Principal Component Analysis (PCA) (see \Cref{sec:princ-comp}) share the commonality of dimensionality reduction, LDA tends to be more robust than PCA since it takes into account the data labels in computing the optimal projection matrix \cite{belhumeur1996eigenfaces}.

Given the dataset $(X,y)$ where each data sample \JSrev{$x_i\in \mathbb{R}^p$} belongs to one of $K$ classes such that if $x_i$ belongs to the $k$-th class then $y_i(k)$ is 1 where $y_i\in \{0,1\}^K$, the input data is partitioned into $K$ groups $\{\pi_k\}^K_{k=1}$ where $\pi_k$ denotes the sample set of the $k$-th class which contains $n_k$ data points. LDA maps the features space \JSrev{$x_i \in \mathbb{R}^p$} to a lower dimensional space $q_i \in \mathbb{R}^r$ ($r< p$) through a linear transformation $q_i=G^\top x_i$ \cite{wang2010multi}. The class mean of the $k$-th class is given by  $\mu_k=\frac{1}{n_k}\sum_{x_i\in \pi_k} x_i$ while the global mean in given by $\mu=\frac{1}{n}\sum_{i=1}^nx_i$. In the projected space the class mean is given by  $\overline{\mu}_k=\frac{1}{n_k}\sum_{q_i\in \pi_k} q_i$ while the global mean in given by $\overline{\mu}=\frac{1}{n}\sum_{i=1}^nq_i$.

The within-class scatter and the between-class scatter 
evaluate the class separability and are defined as $S_w$ and $S_b$ 
respectively such that
\begin{align}
S_w & = \sum_{k=1}^K \sum_{x_i\in \pi_k} (x_i-\mu_k)(x_i-\mu_k)^\top\revtwoCG{,}\\
S_b & = \sum_{k=1}^K n_k (\mu_k-\mu)(\mu_k-\mu)^\top.
\end{align}
The within-class scatter evaluates the spread of the data around the class mean while the between-class scatter evaluates the spread of the class means around the global mean. For the projected data, the within-class and the between-class scatters \JSrev{are defined as $\overline{S}_w$ and $\overline{S}_b$
respectively such that}
\begin{align}
 \overline{S}_w &= \sum_{k=1}^K \sum_{q_i\in \pi_k} (q_i-\overline{\mu}_k)(q_i-\overline{\mu}_k)^\top = G^\top S_wG\revtwoCG{,} \label{within_class}\\
 \overline{S}_b & = \sum_{k=1}^K n_k (\overline{\mu}_k-\overline{\mu})(\overline{\mu}_k-\overline{\mu})^\top = G^\top S_bG. \label{between_class}
\end{align}

The LDA optimization problem is bi-objective where the within-class \revtwoCGblue{scatter} should be minimized while the between-class \revtwoCGblue{scatter} should be maximized. The optimal transformation $G$ can be obtained by maximizing the Fisher criterion (the ratio of between-class to within-class scatters) 
\begin{align}
\max  \frac{|G^TS_bG|}{|G^TS_wG|}. \label{lda_optimization1}
\end{align}
Note that since the between-class and the within-class scatters are not scalar, the determinant is used to obtain a scalar objective function. As discussed in \cite{fukunaga2013introduction}, assuming that $S_w$ is invertible and non-singular, the Fisher criterion is optimized by selecting the $r$ largest eigenvalues of \JSrev{$S_w^{-1}S_b$} and the corresponding eigen vectors \JSrev{$G^*_1,G^*_2,\dots,G^*_r$} form the optimal transformation matrix \JSrev{$G^*=[G^*_1|G^*_2|\dots|G^*_r]$}. \JSrev{Instead of using Fisher criterion, bi-objective optimization techniques may also potentially be used to formulate and solve the LDA optimization problem exactly}.

An alternative formulation of the LDA optimization problem is provided in \cite{chen2013complete} by maximizing the minimum distance between each class center and the total class center. The proposed approach known as the large margin linear discriminant analysis requires the solution of non-convex optimization problems. A solution approach is also proposed based on solving a series of convex quadratic optimization problems. 

%


\subsection{Decision Trees}\label{sec:dec-tress}

Decision trees are classical models for making a decision or classification using splitting rules organized into \revtwoCG{a} tree data structure. Tree-based methods are non-parametric models that partition the predictor space into sub-regions and then yield a prediction based on statistical indicators (e.g., median and mode) of the segmented training data. Decision trees can be used for both regression and classification problems. 

For regression trees, the splitting of the training dataset into distinct and non-overlapping regions can be done using a top-down recursive binary splitting procedure. \revtwoJS{Starting from a root node that contains the full dataset, a cut that splits the data into distinct sets is identified. For the case of a univariate cut (i.e., involving only a single feature), the cutpoint $b$ for feature $j$ is the one that leads to the two splitted regions $R_1=\{x_i| x_{ij} < b\}$ and $R_2=\{x_i| x_{ij} \geq b\}$} that have the greatest possible reduction in the  \mbox{residual sum of squares $\displaystyle\sum_{i: x_i \in R_1(j,b)} (y_i - \hat{y}_{R_1})^2 + \sum_{i: x_i \in R_2(j,b)} (y_i - \hat{y}_{R_2})^2,$} where $\hat{y}_R$ denotes the mean response for the training observations in region $R$. A multivariate split is of the form $a^T x_i < b$, where $a$ is a vector. Another optimization criterion is the measure of purity~\cite{breiman1984classification} such as Gini's index in classification problems. 
For classification problems, \cite{breiman1984classification} highlights that, given their greedy nature, the classical methods based on recursive splitting do not lead to the global optimality of the decision tree. Since building optimal binary decision trees is known to be $\mathcal{NP}$-hard \cite{HYAFIL197615}, heuristic approaches based on mathematical programming paradigms, such as linear optimization~\cite{bennett1992decision}, continuous optimization \cite{bennett1996optimal}, \revtwoJS{and} dynamic programming \cite{6982798, AZAD2014368, cox1989heuristic, 1674938}, have been proposed.

\JSrev{ To find provably optimal decision trees, \cite{bertsimas2017optimal} proposes a mixed integer programming formulation that has an exponential complexity in the depth of the tree. Given a fixed depth $D$, the maximum number of nodes is $T=2^{D+1} - 1$ indexed by $t=1,\dots,T$. Following the notation of \cite{bertsimas2017optimal}, \revtwoJS{the nodes are split} into two sets, branch nodes and leaf nodes. The branch nodes $T_B=\{1,\dots,\lfloor \revtwoJS{\frac{T}{2}} \rfloor\}$ apply a linear split $a^\top x_i < b$ where the left \revtwoJS{child node} includes the \revtwoJS{data points} that satisfy this split while the right \revtwoJS{one} includes the remaining data. In \cite{bertsimas2017optimal}, the splits that are applied at the branch nodes are restricted to a single variable with the option of not splitting a node. \revtwoJS{The binary decision variable $d_t$ takes a value of 1 if branch node $t$ is split and 0 otherwise. Since the splits are univariate, then variable $a_{jt}$, which denotes the value of the coefficient of feature $j$ in the split at node $t$, is also binary. The cutpoint at node $t$ is $b_t \geq 0$.} 

At \revtwoJS{each of} the leaf nodes $\revtwoJS{T_L}=\{\lfloor \revtwoJS{\frac{T}{2}} \rfloor +1, \dots, T\}$, a class prediction is made \revtwoJS{based on} the \revtwoJS{data points} that are included. \revtwoJS{The binary variable $z_{it}$ indicates if data point $i$ is included to leaf node $t$, i.e.\revtwoCGblue{,} $z_{it}=1$ or otherwise $z_{it} = 0$. The binary decision variable $c_{kt}$ takes a value of 1 if leaf node $t$ is assigned label $k$, and 0 otherwise while binary variable $l_t$ indicates if leaf node $t$ is used, i.e.\revtwoCGblue{,} $l_t=1$ or otherwise $l_t = 0$. } 

The mixed integer programming formulation is
 \begin{align}
  \min &\ \frac{1}{\hat{L}}\sum_{t\in T_L}L_t + \alpha \sum_{t\in T_B}d_t \label{dec_tree1}\\
  \mbox{s.t. }  &L_t\geq N_t - N_{kt}-n(1-c_{kt}), \quad \forall k=1,\dots,K, \  t \in T_L\revtwoCG{,} \label{dec_tree2}\\
 &0 \leq L_t\leq N_t - N_{kt}+nc_{kt} \quad \forall k=1,\dots,K, \  t \in T_L\revtwoCG{,} \label{dec_tree3}\\
  &N_{kt}=\frac{1}{2}\sum_{i=1}^n(\revtwoJS{1}+Y_{ik})z_{it}, \quad \forall k=1,\dots,K,\ t\in T_L\revtwoCG{,}\label{dec_tree5}\\
 &N_{t}=\sum_{i=1}^nz_{it} \quad \forall t\in T_L\revtwoCG{,}\label{dec_tree6}\\
 &\sum_{k=1}^Kc_{kt}=l_t \quad \forall t\in T_L\revtwoCG{,}\label{dec_tree7}\\
 &\sum_{t\in T_L} z_{it}=1 \quad \forall i =1,\dots, n\revtwoCG{,}\label{dec_tree10}\\
  &z_{it}\leq l_t \quad \revtwoJS{\forall i =1,\dots, n},\ t\in T_L\revtwoCG{,}\label{dec_tree11}\\
  &\sum_{i=1}^nz_{it}\geq N_{\min}l_t \quad \forall t\in T_L\revtwoCG{,}\label{dec_tree12}\\
  & \revtwoJS{a_m^\top( x_i + \epsilon) \leq b_m + (1+\epsilon_{\max})(1-z_{it})} \quad \forall i=1,\dots,n, \  t\in T_L, \  m\in A_L(t),\label{dec_tree8}\\
 & a_m^\top x_i \geq b_m \revtwoJS{- (1-z_{it})} \quad \forall i=1,\dots,n, \ t\in T_L, \ \forall m\in A_R(t),\label{dec_tree9}\\
   &\sum_{j=1}^pa_{jt}=d_t \quad \forall  t\in T_B\revtwoCG{,}\label{dec_tree13}\\
  &0\leq b_t \leq d_t \quad \forall t\in T_B\revtwoCG{,}\label{dec_tree14}\\
  &d_t\leq d_{p(t)} \quad \forall  t\in T_B\setminus \{1\}\revtwoCG{,}\label{dec_tree15}\\
  & \revtwoJS{z_{it}, \ l_t \in\{0,1\}\quad \forall i=1,\dots,n, \ \forall t \in T_L}\revtwoCG{,}\label{dec_tree16}\\
  & \revtwoJS{c_{kt} \in\{0,1\}\quad \forall k=1,\dots,K, \  t \in T_L,}\\
  & a_{jt}, \revtwoJS{\ d_t} \in\{0,1\}\quad \forall  j=1,\dots,p, \  t\in T_B.\label{dec_tree17}
 \end{align}
\revtwoJS{The objective function \eqref{dec_tree1} minimizes the normalized total misclassification loss $\frac{1}{\hat{L}}\sum_{t\in T_L}L_t$ and the decision tree complexity which is given by $\sum_{t\in T_B}d_t$, the total number of nodes that are split}. $\alpha$ is a tuning parameter and $\hat{L}$ is the baseline loss obtained by predicting the most popular class from the entire dataset. Constraints \eqref{dec_tree2}--\eqref{dec_tree3} set the misclassification loss $L_t$ at leaf node $t$ as $L_t = N_t - N_{kt}$ if node $t$ is assigned label $k$ (i.e $c_{kt} = 1$), where $N_t$ is the total number of data points at leaf node $t$ and $N_{kt}$ is the total number of data points at node $t$ whose true labels are~$k$. The counting of $N_{kt}$ and $N_t$ is enforced by \eqref{dec_tree5} and \eqref{dec_tree6}, respectively, \revtwoJS{where $Y_{ik}$ is a parameter taking the value of 1 if data point $i$ has a label $k$ and $-1$ otherwise}. Constraints \eqref{dec_tree7} indicate that each leaf node that is used (i.e.\revtwoCGblue{,} $l_t=1$) should be assigned to a label $k = 1\dots K$. Constraints \eqref{dec_tree10} indicate that each data point should be assigned to exactly one leaf node. 
Constraints \eqref{dec_tree11}--\eqref{dec_tree12} indicate that data points can be assigned to a node only if that node is used and if a node is used then at least $N_{\min}$ data points should be assigned to it. The splitting of the data points at each of the branch nodes is enforced by constraints \eqref{dec_tree8}--\eqref{dec_tree9}
where $A_L(t)$ is the set of ancestors of $t$ whose left branch has been followed on the path from the root node to node~$t$. Similarly, $A_R(t)$ is the set of ancestors of $t$ whose right branch has been followed on the path from the root node to node~$t$. $\epsilon$ and $\epsilon_{\max}$ are small numbers to enforce the strict split $a^\top x < b$ at the left branch (see \cite{bertsimas2017optimal} for finding good values for $\epsilon$ and $\epsilon_{\max}$). Constraints \eqref{dec_tree13}--\eqref{dec_tree14} indicate that the splits are restricted to a single variable with the option of not splitting a node ($d_t=0$). As enforced by constraints \eqref{dec_tree15}, if $p(t)$, the parent of node $t$, does not apply a split then so is node $t$. Finally constraints \eqref{dec_tree16}--\eqref{dec_tree17} set the binary conditions.}

 \ \\
 
	 An alternative formulation to the optimal decision tree problem is provided in \cite{gunluk2018optimal}. The main difference between the formulation of \cite{gunluk2018optimal} and \cite{bertsimas2017optimal} is that the approach of \cite{gunluk2018optimal} is specialized to the case where the features take categorical values. By exploiting the combinatorial structure that is present in the case of categorical variables, \cite{gunluk2018optimal} provides a strong formulation of the optimal decision tree problem thus improving the computational performance. Furthermore the formulation of \cite{gunluk2018optimal} is restricted to binary classification and the tree topology is fixed, which lowers the required computational effort for solving the optimization problem to optimality. \JSrev{A commonality between the models presented in \cite{bertsimas2017optimal} and \cite{gunluk2018optimal} is that the split that is considered at each node of the decision tree involves only one variable mainly to achieve better computational performance when solving the optimization model. More generally, splits that span multiple variables can also be considered at each node as presented in \cite{blanquero2018optimal, verwer2017learning,verwer2017auction}. The approach of \cite{blanquero2018optimal}\revtwoCG{,} which is extended in \cite{blanquero2018sparsity} to account for sparsity by using regularization, is based on a nonlinear continuous optimization formulation to learn decision trees with general splits.}


While single decision tree models are often preferred by data analysts for their high interpretability, the model accuracy can be largely improved by taking multiple decision trees into account. Such approaches include bagging, random forests, and boosting. \JSrev{Bagging creates multiple decision trees by obtaining several training subsets by randomly choosing with replacement data points from the training set and subsequently training a decision tree for each subset. Random forests \revtwoCGblue{create} training subsets similar to bagging with the addition of randomly selecting a subset of features for training each tree. Boosting iteratively creates decision trees where a weight on the training data is set and is increased at each iteration for the misclassified data points so as to subsequently create a decision tree that is more likely to correctly classify previously misclassified data. \revtwoCGblue{These} models that make predictions based on aggregating the predictions of individual trees are also known as tree ensemble. A mixed integer optimization model for tree ensemble has been recently proposed in \cite{mivsic2017optimization}.}

Decision trees can also be used in a more general range of applications as algorithms for problem solving, data mining, and knowledge representation. In \cite{7321423}, several greedy and dynamic programming approaches are compared for building decision trees on datasets with inconsistent labels (i.e, many-valued decision approach). Many-valued decisions can be evaluated in terms of multiple cost functions in a multi-stage optimization \cite{AZAD2017910}. Recently, \cite{CHIKALOV2018689} investigated conflicting objectives in the construction of decision trees by means of bi-criteria optimization. Since the single objectives, such as minimizing average depth or the number of terminal nodes, are known to be $\mathcal{NP}$-hard, the authors propose a bi-criteria optimization approach by means of dynamic programming.

\subsection{Support Vector Machines}\label{sec:svm}

Support vector machines (SVMs) are another class of supervised machine learning algorithms that are based on statistical learning and \revtwoCG{have} received significant attention in the optimization literature \cite{carrizosa2013supervised, vapnik1998statistical, vapnik2013nature}. \JSrev{Given a training set $(X,y)$ with $n$ training inputs where $X\in \mathbb{R}^{n\times p}$ and binary response variables $y\in \{-1,1\}^n$,} the objective of the support vector machine problem is to identify a hyperplane \revtwoJS{$w^\top x + \gamma = 0$, where $w\in \mathbb{R}^p$ and $\gamma\in \mathbb{R}$, which} separates the two classes of data points with a maximal separation margin measured as  the width of the band that separates the two classes. \revtwoJS{In this section, $w$ denotes the vector of coefficients corresponding to each of the input variables and $\gamma$ is the intercept of the separating hyperplane.} As detailed next, the underlying optimization problem is a \JSrev{linearly constrained} convex quadratic optimization problem.

\subsubsection{Hard Margin SVM}
The most basic version of SVMs is the hard margin SVM that assumes that there exists a hyperplane that geometrically separates the data points into the two classes such that no data point is misclassified \cite{Cortes1995}. The training of the SVM model involves finding the hyperplane  that separates the data and whose distance to the closest data point in either of the classes, i.e.\revtwoCGblue{,} margin, is maximized. 

The distance of a particular data point $x_i$ to the separating hyperplane is
\begin{align*}
 \frac{y_i(w^\top x_i + \gamma)}{\|w\|_2}\revtwoCG{.}
\end{align*}
\revtwoJS{The distance to the closest data point is normalized to $\frac{1}{\|w\|_2}$ where $\|w\|_2$ denotes the $2$-norm}. Thus the data points with labels $y=-1$ are on one side of the hyperplane such that $w^\top x + \gamma \leq 1$ while the data point with labels $y=1$ are on the other side $w^\top x + \gamma \geq 1$. The optimization problem for finding the separating hyperplane is then
\begin{align*}
 \max &  \ \frac{1}{\|w\|_2}\\
 \mbox{s.t. }& y_i(w^\top x_i + \gamma) \geq 1 \quad \forall  i = 1,\dots,n\revtwoCG{,}\\
 & \revtwoJS{w \in \mathbb{R}^p, \gamma \in \mathbb{R},} 
\end{align*}
which is equivalent to
\begin{align}
\min &  \ \|w\|_2^2\label{hardsvm_1}\\
 \mbox{s.t. }& y_i(w^\top x_i + \gamma) \geq 1 \quad \forall  i = 1,\dots,n\revtwoCG{,}\label{hardsvm_2}\\
 & \revtwoJS{w \in \mathbb{R}^p, \gamma \in \mathbb{R},} \label{hardsvm_3}
\end{align}
that is a convex quadratic problem.

Forcing the data to be separable by a linear hyperplane is a strong condition that often does not \revtwoCG{hold in practice. Thus, the soft-margin SVM, which relaxes the condition of perfect separability, is used instead.}

\subsubsection{Soft-Margin SVM}
When the data is not linearly \revtwoJS{separable}, problem \eqref{hardsvm_1}--\eqref{hardsvm_3} is infeasible. \revtwoJS{Alternatively, \cite{bennett1992robust} proposed a linear program that minimizes a weighted average of the errors given by the points lying on the wrong  side of the separator.} \JSrev{This work was then extended in \cite{Cortes1995} which  presented the soft margin SVM. The soft margin SVM introduces slack variables \revtwoJS{$\xi_i \geq 0$} into constraints \eqref{hardsvm_2} which are then} penalized in the objective function as a proxy to minimizing the number of data points that are on the wrong side. The soft-margin SVM optimization problem is
\begin{align}
\min &  \ \|w\|_2^2+C\sum_{i=1}^n\xi_i \label{softsvm_1}\\
 \mbox{s.t. }& y_i(w^\top x_i + \gamma) \geq 1 - \xi_i \quad \forall i = 1,\dots,n\revtwoCG{,} \label{softsvm_2}\\
 & \revtwoJS{w \in \mathbb{R}^p, \gamma \in \mathbb{R},}\\
 &\xi_i \geq 0\quad \forall i = 1,\dots,n.\label{softsvm_3}
\end{align}

Another common alternative is to include the error term $\xi_i$ in the objective function by using the squared hinge loss $\sum_{i}^n\xi_i^2$ instead of the hinge loss $\sum_{i}^n\xi_i$. \JSrev{The hinge loss function takes a value of zero for a data point that is correctly classified while it takes a positive value that is proportional to the distance from the separating hyperplane for a misclassified data point.} Hyperparameter $C$ is then tuned to obtain the best classifier.

Besides the direct solution of problem \eqref{softsvm_1}--\eqref{softsvm_3} as a convex quadratic problem, replacing the \CGrev{$2$-norm by the $1$-norm} leads to a linear optimization problem generally at the expense of higher misclassification rate \cite{bradley2000massive}. 

\subsubsection{Sparse SVM}
Using the \CGrev{$1$-norm} is also an approach to sparsify $w$, i.e.\revtwoCGblue{,} reduce the number of features that are involved in the classification model \cite{bradley2000massive, zhu20041}. An approach known as the elastic net includes both the \CGrev{$1$-norm and the $2$-norm} in the objective function and tunes the bias towards one of the norms through a hyperparameter \cite{wang2006doubly, zou2005regularization}. \JSrev{Several other approaches for dealing with sparsity in SVM have been proposed in \cite{aytug2015feature, DUNBAR2010470, gaudioso2017lagrangian, ghaddar2018high, guyon2002gene, maldonado2014feature,PIRAMUTHU2004483}.} \JSrev{The number of features can be explicitly modeled in \eqref{softsvm_1}--\eqref{softsvm_3} by using binary variables $z\in\{0,1\}^p$ where $z_j = 1$ indicates that feature $j$ is selected and otherwise $z_j = 0$ \cite{chan2007direct}. A constraint limiting the number of features to a maximum desired number can be enforced resulting in the following mixed integer quadratic problem
\begin{align}
\min &  \ \|w\|_2^2+C\sum_{i=1}^n\xi_i \label{sparsesvm_1}\\
 \mbox{s.t. }& y_i(w^\top x_i + \gamma) \geq 1 - \xi_i \quad \forall i = 1,\dots,n\revtwoCG{,} \label{sparsesvm_2}\\
 & -Mz_j\leq w_j \leq Mz_j \quad \forall j=1,\dots,p\revtwoCG{,}\label{sparsesvm_3}\\
 &\sum_{j=1}^pz_j \leq r\revtwoCG{,}\label{sparsesvm_4}\\
& \revtwoJS{w \in \mathbb{R}^p, \gamma \in \mathbb{R},} \\
 &z_j\in\{0,1\}\quad \forall j = 1,\dots,p\revtwoCG{,}\label{sparsesvm_5}\\
 &\xi_i \geq 0\quad \forall i = 1,\dots,n.\label{sparsesvm_6}
\end{align}
Constraints \eqref{sparsesvm_3} force $z_j = 1$ when feature $j$ is used, i.e.\revtwoCGblue{,} $w_j \neq 0$ ($M$ denotes a sufficiently large number). Constraints \eqref{sparsesvm_4} set a limit $r$ on the maximum number of features that can be used.}

\subsubsection{The Dual Problem and Kernel Tricks}
The data points can be mapped to a higher dimensional space through a mapping function $\phi(x)$. A soft margin SVM can then be applied such that
\begin{align}
\min &  \ \|w\|_2^2+C\sum_{i=1}^n\xi_i\label{mapped_svm1}\\
 \mbox{s.t. }& y_i(w^\top \phi(x_i) + \gamma) \geq 1 - \xi_i \quad \forall i = 1,\dots,n\revtwoCG{,}\label{mapped_svm2}\\
 & \revtwoJS{w \in \mathbb{R}^p, \gamma \in \mathbb{R},}\\
 &\xi_i \geq 0\quad \forall i = 1,\dots,n.\label{mapped_svm3}
\end{align}
Through this mapping, the data has a linear classifier in the higher dimensional space however a \CGrev{nonlinear} separation function is obtained in the original space.

To solve problem \eqref{mapped_svm1}--\eqref{mapped_svm3}, the following dual problem is first obtained
\begin{align*}
 \max_\alpha \ & \sum_{i=1}^n\alpha_i-\frac{1}{2}\sum_{i,j=1}^n\alpha_i \alpha_jy_iy_j \phi(x_i)^\top \phi(x_j)\\
  \mbox{s.t. }& \sum_{i=1}^n\alpha_iy_i = 0 \quad \forall i= 1,\dots,n\revtwoCG{,}\\
 & 0\leq \alpha_i \leq C \quad \forall i= 1,\dots,n\revtwoCG{,}
\end{align*}
where $\alpha_i$ are the dual variables of constraints \eqref{mapped_svm2}. Given a kernel function $K:\mathbb{R}^m\times\mathbb{R}^m\rightarrow \mathbb{R}$ where $K(x_i,x_j)=\phi(x_i)^\top \phi(x_j)$, the dual problem is 
\begin{align*}
 \max_\alpha \ & \sum_{i=1}^n\alpha_i-\frac{1}{2}\sum_{i,j=1}^n\alpha_i \alpha_jy_iy_j K(x_i,x_j)\\
  \mbox{s.t. }& \sum_{i=1}^n\alpha_iy_i = 0,\quad \forall i= 1,\dots,n\revtwoCG{,}\\
 & 0\leq \alpha_i \leq C, \quad \forall i= 1,\dots,n\revtwoCG{,}
\end{align*}
which is a convex quadratic optimization problem. Thus only the kernel function $K(x_i,x_j)$ is required while the explicit mapping $\phi$ is not needed. 

\revtwoJS{Among the commonly used kernel functions is the polynomial $K(x_i,x_j)=(x_i^\top x_j +c)^d$ where $c$ controls the trade-off between the higher-order and the lower-order terms in the polynomial and $d$ is the degree of the polynomial. The polynomial kernel models the interaction between the data up to the degree $d$. A high degree polynomial tends to overfit the training data. Another kernel function is the radial basis $K(x_i,x_j)=e^{-\frac{\|x_i-x_j\|_2^2}{\gamma}} $ where $\gamma$ acts as a smoothing parameter. A smaller $\gamma$ tends to overfit the training data. The sigmoidal kernel $K(x_i,x_j)=\tanh(\varphi {x_i}^\top x_j+c)$ is also commonly used where $\varphi$ is a scaling parameter of the input data and $c$ is a shifting parameter that controls the threshold of the mapping. Further details on kernel functions are provided in  \cite{alam2018kernel,carrizosa2013supervised, herbrich2001learning}.}

\cg{Since the classification in high dimensional space can be difficult to interpret for practitioners, Binarized SVM (BSVM) replaces the continuous predictor variables with a linear combination of binary cutoff variables \cite{CarrIJOC}. \revtwoJS{BSVM is also extended in \cite{CARRIZOSA2011260} to capture the interactions between the relevant variables}.} \js{Another important practical aspect to consider is data uncertainty. Often the training data suffers from inaccuracies in the labels and in the features that are collected which may negatively affect the performance of the classifiers. While typically regularization is used to mitigate the effect of uncertainty, \cite{bertismas2019robustclassification} proposes robust optimization models for logistic regression, decision trees, and support vector machines and shows increased accuracy over regularization\revtwoCGblue{, and} most importantly without changing the complexity of the classification problem.}

\JSrev{\subsubsection{Support Vector Regression}
Although as discussed earlier, \revtwoCG{SVM has} been introduced for binary classification, its extension to regression, i.e.\revtwoCGblue{,} support vector regression, has received significant interest in the literature \cite{smola2004tutorial}. The core idea of support vector regression is to find a linear function $f(x)=w^\top x + \gamma$ that can approximate with a tolerance $\epsilon$ a training set $(X,y)$ where $y\in \mathbb{R}$ \cite{vapnik2013nature}. Such a linear function may however not exist, and thus \revtwoJS{slack variables $\xi^+_i \geq 0$ and $\xi^-_i \geq 0$ denoting positive and negative deviations from the desired tolerance} are introduced and minimized similar to the soft-margin SVM. The corresponding optimization problem is
\begin{align}
 \min\ & \|w\|_2^2+C\sum_{i=1}^n(\xi^+_i + \xi^-_i)\label{softsvr_1}\\
 \mbox{s.t. }& y_i - w^\top x_i -\gamma \leq \epsilon + \xi^+ \quad \forall i = 1,\dots,n\revtwoCG{,}\\
 & w^\top x_i +\gamma - y_i\leq \epsilon + \xi^- \quad \forall i = 1,\dots,n\revtwoCG{,}\\\
 & \revtwoJS{w \in \mathbb{R}^p, \gamma \in \mathbb{R},}\\
 &\xi^+_i,\ \xi^-_i \geq 0\quad \forall i = 1,\dots,n.\label{softsvr_4}
\end{align}
Hyperparameter $C$ is tuned to adjust the weight on the deviation from the tolerance $\epsilon$. This deviation from $\epsilon$ is the $\epsilon$-insensitive loss function $|\xi|_\epsilon$ given by
$$|\xi|_\epsilon=\begin{cases}
         0 & \mbox{ if }|\xi|\leq \epsilon\revtwoCG{,}\\
         |\xi|-\epsilon &\mbox{ otherwise.}
        \end{cases}
$$
As detailed extensively in \cite{smola2004tutorial}, kernel tricks can also be applied to \eqref{softsvr_1}--\eqref{softsvr_4} which is solved by formulating the dual problem.
}

\subsubsection{Support Vector Ordinal Regression}

\JSrev{In situations where the data contains ordering preferences, i.e.\revtwoCGblue{,} the training data is labeled \revtwoCG{by ranks, where} the order of the rankings is relevant while the distances between the ranks is not defined or irrelevant to the training, the purpose of learning is to find a model that maps the preference information.}

\JSrev{The application of classic regression models for such type of data requires the transformation of the ordinal ranks to \revtwoCG{numerical values. However,} such approaches often fail in providing robust models as an appropriate function to map the ranks to distances is challenging to find \cite{kramer2001}. An alternative is to encode the ordinal ranks into binary classifications at the expense of a large increase in the scale of the problems \cite{har2003constraint, herbrich2000}. 

An extension of SVM for ordinary data has been proposed in \cite{shashua2003ranking} and extended in \cite{chu2007support}. Given a training dataset with $r$ ordered categories $\{1,\dots,r\}$ where $n_j$ is the number of data points labeled as order $j$, the support vector ordinal regression finds $r-1$ separating parallel hyperplanes $w^{\top}x+\beta_j = 0$ where $\beta_j$ is the threshold associated with the hyperplane that separates the orders $k\leq j$ from the remaining orders. Thus $x_{i,k}$, the $i^{\mbox{th}}$ data sample of order $k \leq j$, should have a function value lower than the margin $\beta_j-1$ while the data samples with orders $k > j$ should have a function value greater than the margin $\beta_j+1$. The errors for violating these conditions are given by $\xi^+_{i,kj}\geq 0$ and $\xi^-_{i,kj}\geq 0$ respectively. Following \cite{chu2007support}, the associated SVM formulation is
\begin{align*}
 \min &  \ \|w\|_2^2+C\sum_{j=1}^{r-1}(\sum_{k=1}^j\sum_{i=1}^{n_k}\xi^+_{i,kj} + \sum_{k=j+1}^r\sum_{i=1}^{n_k}\xi^-_{i,kj})\\
 \mbox{s.t. } & w^\top x_{i,k} - \beta_j \leq -1 + \xi^+_{i,kj} \quad \forall k = 1, \dots, j, \  j= 1,\dots, r-1 \revtwoCG{,} \  i=1,\dots,\revtwoJS{n_k}\revtwoCG{,}\\
    & w^\top x_{i,k} - \beta_j \geq 1 - \xi^-_{i,kj} \quad \forall k = j+1, \dots, r,\  j= 1,\dots,r-1, \  i=1,\dots,\revtwoJS{n_k},\\
    & \revtwoJS{w \in \mathbb{R}^p, \ \beta_j \in \mathbb{R} \quad \forall j= 1,\dots,r-1,}\\
    & \revtwoJS{\xi^+_{i,kj} \geq 0 \quad \forall k = 1, \dots, j, \  j= 1,\dots, r-1, \  i=1,\dots,\revtwoJS{n_k},}\\
    & \revtwoJS{\xi^-_{i,kj} \geq 0 \quad \forall k = j+1, \dots, r,\  j= 1,\dots,r-1, \  i=1,\dots,n_k.}
\end{align*}

As detailed in \cite{chu2007support}, kernel tricks can be also applied by considering the dual problem. Finally we note that preference modeling using machine learning has several commonalities with various approaches in multi-criteria decision analysis and most notably, robust ordinal regression. We refer the readers to \cite{corrente2013robust} for a detailed comparison between preference learning using machine learning and muti-criteria decision making.}

\section{Clustering}\label{sec:clustering}

Data clustering is a class of unsupervised learning approaches that has been widely used, particularly in applications of data mining, pattern recognition, and information retrieval. \revtwoJS{Given an input $X\in \mathbb{R}^{n\times p}$, which includes $n$ unlabeled observations $x_1, \dots, x_n$ with $x_i \in \mathbb{R}^{p}$,} cluster analysis aims at finding  $K$ subsets 
of $X$, called clusters, which are homogeneous and well separated. \textit{Homogeneity} indicates the similarity of the observations within the same cluster (typically, by means of a distance metric), while the \textit{separability} accounts for the differences between entities of different clusters. The two concepts can be measured via several criteria and lead to different types of clustering algorithms (see, e.g., \cite{hansen1997cluster}). The number of clusters is typically a tuning parameter to be fixed before determining the clusters.
An extensive survey on data clustering analysis is provided in \cite{jain1999data}.

In case the entities are points in a Euclidean space, the clustering problem is often modeled as a  network problem and shares many similarities with classical problems in operations research, such as the $p$-median problem \CGrev{\cite{benati2014mixed, klastorin1985p, MAI2018594, mulvey1979cluster}. 
}
 In the following subsections, the commonly used \cg{minimum sum-of-squares clustering, the capacitated clustering, and the $K$-hyperplane clustering are discussed.}

\subsection{Minimum Sum-Of-Squares Clustering (a.k.a. $K$-Means Clustering)}\label{sec:sos-clustering}
Minimum sum-of-squares clustering is one of the most commonly adopted clustering algorithms. It requires to find a number of disjoint clusters for observations \CGrev{$x_i, i=1, \dots, n$, where $x_i\in\mathbb{R}^p$} such that the distance to cluster centroids is minimized. Given that typically the number of clusters $K$ is a-priori fixed, the problem is also referred to as $K$-means clustering. \cg{The decision of the cluster size is typically taken by examining the elbow curve, or similarity indicators, such as silhouette values and Calinski-Harabasz index, or via mathematical programming approaches including the maximization of the modularity of the associated graph \cite{cafieri2014reformulation, CAFIERI201465}.} 

Defining the binary variables 
$$u_{ij}=\begin{cases}
         1 &\mbox{ if observation $i$ belongs to cluster $j$}\\
         0 &\mbox{ otherwise,}
        \end{cases}
$$
\JSrev{and the centroid $\mu_j \in \mathbb{R}^p$ of each cluster $j$,} the problem \CGrev{of minimizing the within-cluster variance} is formulated in \cite{Aloise2012} as the following \CGrev{mixed integer} nonlinear program 
\begin{align}
\min \ &  \sum_{i = 1}^n\sum_{j=1}^K u_{ij} \lVert x_i - \mu_j \rVert_2^2   \label{k-means:obj}\\
\mbox{s.t. } & \sum_{j =1}^K u_{ij} = 1 \quad \forall i = 1,\dots,n\revtwoCG{,}\label{k-means:con1}\\
& \mu_{j} \in \mathbb{R}^p \quad \forall j = 1,\dots,K\revtwoCG{,}\label{k-means:con2}\\
& u_{ij} \in \{0,1\} \quad \forall i = 1,\dots,n, \ j = 1,\dots,K.\label{k-means:con3}
\end{align}
\JSrev{By introducing the variables $d_{ij}$ which denote the distance of observation $i$ from centroid $j$}, the following linearized formulation is obtained 
\begin{align*}
 \min \ & \sum_{i = 1}^n\sum_{j=1}^K d_{ij}\\
 \mbox{s.t. } & \sum_{j =1}^K u_{ij} = 1, \quad \forall i = 1,\dots,n\revtwoCG{,}\\
  & d_{ij} \geq ||x_i-\mu_j||_2^2 - M(1-u_{ij})  \quad \forall i=1,\dots,n, \ j = 1,\dots,K\revtwoCG{,}\\
  & \mu_{j} \in \mathbb{R}^p \quad \forall j =1,\dots,K\revtwoCG{,}\\
  & u_{ij} \in \{0,1\}, \ d_{ij}\geq 0 \quad \forall i=1,\dots,n, \ j = 1,\dots,K.
\end{align*}
\JSrev{Parameter $M$ is a sufficiently large number. A \revtwoJS{heuristic} solution approach based on the gradient method is proposed for problem \eqref{k-means:obj}--\eqref{k-means:con3} in \cite{bagirov2006new}. Alternatively, a column generation approach for large-scale instances has been proposed in \cite{Aloise2012} and a bundle approach has been presented in \cite{karmitsa2017new}}.

The case where the space is not Euclidean is considered in \cite{CARRIZOSA2013356}. Alternatively, \cite{santi2016model} presents the Heterogeneous Clustering Problem (HCP) where the observations to cluster are associated with multiple dissimilarity matrices\revtwoCG{.} HCP is formulated as a mixed integer quadratically constrained quadratic program.
Another variant is presented in \cite{sauglam2006mixed} where the homogeneity is expressed by the minimization of the maximum diameter $D_{\text{max}}$ of the clusters\cg{. The resulting nonconvex bilinear mixed integer program is solved via a graph-theoretic approach based on seed finding.}

\CGrev{
Many common solution approaches for $K$-means clustering are based on heuristics. A popular method implemented in data science packages (e.g., scikit-learn \cite{scikit-learn}) is the two-step improvement procedure proposed in \cite{macqueen1967some}. Starting from a sample of $K$ points in set $X$ as initial cluster centers (centroids $\mu_j^0$), at each iteration $k$, the algorithm assigns each point in $X$ to the nearest centroid $\mu_j^k$ and then computes the centroids $\mu_j^{k+1}$ of the new partition. The procedure is guaranteed to decrease the 
within-cluster variance and it is run until this metric is sufficiently low. Given the dependency of the procedure to the choice of $\mu_j^0$, typically the clustering is repeated with different initial centroids and the best clusters are selected. Other heuristics relax the assumption to produce exactly $K$ clusters. For instance, \cite{macqueen1967some} merges clusters if their centroids are sufficiently close. Clustering is also used within heuristics for hard combinatorial problems \revtwoJS{(\cite{ganesh2007cloves, kwatera1993clustering})}, and can be integrated in problems where the evaluation of multiple solutions is important (e.g.\revtwoCGblue{,} Cluster Newton Method \cite{aoki2014cluster, gaudreau2015improvements}). Cluster Newton method approximates the Jacobian in the domain covered by the cluster of points, instead of locally as done by the  traditional Newton's Method~\cite{kelley1999iterative}, and this has a regularization effect.
}

\subsection{Capacitated Clustering}\label{sec:cap-clustering}
The Capacitated Centered Clustering Problem (CCCP) deals with finding a set of clusters with a capacity limitation and homogeneity expressed by the similarity to the cluster centre. Given a set of potential clusters \CGrev{${1,\dots,K}$}, a mathematical formulation for CCCP is given in \cite{negreiros2006capacitated} as 
\begin{align}
\min  \  & \sum_{i=1}^{n} \sum_{j=1}^{K} \underline{s}_{ij} u_{ij}\label{cluster_CCP:obj} \\
\text{s.t } & \sum_{j=1}^{K} u_{ij} = 1 \quad \forall i = 1, \dots, n,\label{cluster_CCP:con1}  \\
& \sum_{j=1}^{K} v_j \leq K \revtwoCG{,} \label{cluster_CCP:con2} \\
& u_{ij} \leq v_j \quad \forall i = 1, \dots, n, \ j = 1, \dots, V\revtwoCG{,} \label{cluster_CCP:con3} \\
& \sum_{i=1}^{n} q_i u_{ij} \leq Q_j \quad \forall j = 1, \dots, K \revtwoCG{,} \label{cluster_CCP:con4} \\
& u_{ij}, v_j \in \{0,1\} \quad \forall i = 1, \dots, n, \ j = 1, \dots, K. \notag
\end{align}
\JSrev{Parameter \revtwoCG{$K$} is an upper bound on the number of clusters, $\underline{s}_{ij}$ is the dissimilarity measure between observation $i$ and cluster $j$, $q_i$ is the weight of observation $i$, and $Q_j$ is the capacity of cluster $j$. Variable $u_{ij}$ denotes the assignment of observation $i$ to cluster $j$ and variable $v_j$ is equal to $1$ if cluster $j$ is used. If the metric $\underline{s}_{ij}$ is a distance and the clusters are homogeneous (i.e., $Q_j = Q \; \forall j$), the formulation also models the well-known facility location problem. A solution approach is discussed in \cite{chaves2010clustering} while an alternative quadratic programming formulation is presented in~\cite{Lewis2014}. Solution heuristics have also been proposed in \cite{MAI2018594} and \cite{SCHEUERER2006533}.}

\subsection{$K$-Hyperplane Clustering}\label{sec:k-hyp-clust}

In the $K$-Hyperplane Clustering ($K$-HC) problem, a hyperplane, instead of a center, is associated with each cluster. This is motivated by applications such as text mining and image segmentation, where collinearity and coplanarity relations among the observations are the main interest of the unsupervised learning task, rather than the similarity. Given the observations \CGrev{$x_i, i=1, \dots, n$}, the $K$-HC problem requires to find $K$ clusters, and a hyperplane \CGrev{$H_j = \{x \in \mathbb{R}^p: w_j^T x = \gamma_j\}$}, with $w_j \in \mathbb{R}^p$ and $\gamma_j \in \mathbb{R}$, \revtwoCG{for each cluster $j$. The aim is to minimize} the sum of the squared $2$-norm Euclidean orthogonal distances between each observation and the corresponding cluster. 

Given that the orthogonal distance of $x_i$ to hyperplane $H_j$ is given by $\frac{|w_j^T x_i - \gamma_j|}{\lVert w \rVert_2}$, $K$-HC is formulated in \cite{amaldi2013distance} as the following \CGrev{mixed integer} quadratically \revtwoCG{constrained} quadratic problem\newpage
\JSrev{\begin{align}
\min  \  & \sum_{i=1}^{n} \delta_i^2   \label{hyperplane_clustering:obj} \\
\text{s.t } & \sum_{j=1}^{K} u_{ij} = 1 \quad \forall i = 1, \dots, n\revtwoCG{,} \label{hyperplane_clustering:con1}  \\
& \delta_i \geq (w_j ^T x_i - \gamma_j) - M (1-u_{ij}) \quad \forall i = 1, \dots, n, \ j = 1, \dots, K\revtwoCG{,} \label{hyperplane_clustering:bigM1} \\
& \delta_i \geq (-w_j ^T x_i + \gamma_j) - M (1-u_{ij}) \quad \forall i = 1, \dots, n, \ j = 1, \dots, K\revtwoCG{,} \label{hyperplane_clustering:bigM2} \\
& \lVert w_j \rVert_2 \geq 1 \quad \forall j = 1, \dots, K\revtwoCG{,} \label{hyperplane_clustering:normw}\\
& \delta_i \geq 0 \quad \forall i = 1, \dots, n\revtwoCG{,} \label{hyperplane_clustering:y}\\
& w_j \in \mathbb{R}^p, \gamma_j \in \mathbb{R} \quad \forall j = 1, \dots, K\revtwoCG{,} \label{hyperplane_clustering:wgamma}\\
& u_{ij} \in \{0,1\}  \quad \forall i \in 1, \dots, n, \ j = 1, \dots, K. \label{hyperplane_clustering:x}
\end{align}
\revtwoCG{Binary variable $u_{ij}$ is equal to $1$ if point $x_i$ is assigned to cluster $j$, and $0$ otherwise. Linear constraints \eqref{hyperplane_clustering:bigM1}--\eqref{hyperplane_clustering:bigM2} set $\delta_i$ as the distance between point $x_i$ and the hyperplane of cluster $j$. These constraints are enforced only if $u_{ij}$ is equal to $1$, otherwise they are redundant.} 
The non-convexity is due to constraint \eqref{hyperplane_clustering:normw}. As a solution approach,  a distance-based reassignment heuristic that outperforms spatial branch-and-bound solvers is proposed in \cite{amaldi2013distance}.}

\section{\CGrev{Linear Dimension Reduction}}\label{sec:dim-red}
\CGrev{In \Cref{sec:shrink}, shrinkage methods have been discussed as a way to improve model interpretability by fitting a model with all original $p$ predictors. In this section, we discuss dimension reduction methods that search for \JSrev{$H<p$} linear combinations of the predictors such that \JSrev{$Z_h = \sum_{j=1}^p \phi_{j}^h X_j$} (also called \textit{projections}) where $X_j$ denotes column $j$ of X, i.e.\revtwoCGblue{,} the vector of values of feature $j$ of the training set.} \JSrev{While this section focuses on Principle Component Analysis and Partial Least Squares, we note that other linear and nonlinear dimension reduction methods exist. An extensive survey  on benefits and shortcomings of dimension reduction methods is presented in \cite{cunningham2015linear}.}

\subsection{Principal Components}\label{sec:princ-comp}

Principal Components Analysis (PCA) \cite{jolliffe2011principal} \revtwoCG{aims to find a low-dimensional representation of the dataset with highly informative derived features.
Principal components are ordered in terms of their explained variances, which measure the amount of information retained from the original set of features $X_1,\ldots,X_p$.}

 In particular, assuming the regressors are standardized to a mean of $0$ and a variance of $1$, the direction of the first principal component is a unit vector \JSrev{$\phi^1 \in \mathbb{R}^p$} that is the solution of the optimization problem
\JSrev{\begin{align}
	\max_{\phi^1 \revtwoCG{\in \mathbb{R}^p}} & \ \frac{1}{n}\sum_{i=1}^n\left(\sum_{j=1}^p\phi_{j}^1x_{ij}\right)^2 \label{PCA1-obj}\\
	\mbox{s.t. }& \sum_{j=1}^p (\phi_{j}^1)^2=1. \label{PCA1-cons}
	\end{align}}
Problem \eqref{PCA1-obj}--\eqref{PCA1-cons} is the traditional formulation of PCA and can be solved via Lagrange multipliers methods. Since the formulation is sensitive to the presence of outliers, several approaches have been proposed to improve robustness \cite{reris2015principal}. One approach is to replace the \revtwoCG{$2$-norm} in \eqref{PCA1-obj} with the \revtwoCG{$1$-norm}.

\JSrev{An iterative approach can be used to obtain the \revtwoCG{principal} components where} the first principal component  \JSrev{$Z_1 = \sum_{j=1}^p \phi_{j}^1 X_j$} is the projection of the original features with the largest variability. \JSrev{The subsequent principal components are obtained iteratively where each principal component \mbox{$Z_h, h = 2, \dots, H$} is obtained by a linear combination of the feature columns $X_1, \dots, X_p$. Each $Z_h$ is uncorrelated with $Z_1, \dots, Z_{h-1}$ which have larger variance}. Introducing the sample covariance matrix $S$ of the regressors $X_j$, the direction \JSrev{$\phi^h \in \mathbb{R}^p$ of the $h$-th principal component $Z_h$} is the solution of 
\newpage
\JSrev{\begin{align}
	\max_{\phi^h \revtwoCG{\in\mathbb{R}^p}} & \ \frac{1}{n}\sum_{i=1}^n\left(\sum_{j=1}^p\phi_{j}^hx_{ij}\right)^2 \label{PCAm-obj}\\
	\mbox{s.t. }& \sum_{j=1}^p (\phi_{j}^h)^2=1 \revtwoCG{,}\label{PCAm-cons-var}\\
	& {\phi^h}^\top S \phi^l = 0 \quad \forall l=1, \dots, h-1. \label{PCAm-cons:uncorr}
	\end{align}}

\CGrev{
	PCA can be used for several data analysis problems which benefit from reducing the problem dimension. Principal Components Regression (PCR) is a two-stage procedure that uses the first principal components as predictors for a linear regression model. PCR has the advantage of including less predictors than the original set and at the same time retaining the variability of the dataset in the derived features.
	However, principal components might not be relevant with the response variables of the regression. 
	
	To select principal components in regression models, the regression loss function and the PCA objective function can be combined in a single-step quadratic programming formulation~\cite{kawano2015sparse}.
}
Since the identification of the principal components does not require any knowledge of the response $y$, PCA can \revtwoCGblue{also} be adopted in unsupervised learning such as in the $k$-means clustering method (see Section \ref{sec:sos-clustering}, \cite{ding2004k}). 
A known drawback of PCA is interpretability. To promote the sparsity of the projected components\CGrev{, and thus make them more interpretable,} \cite{CARRIZOSA2014349} \CGrev{formulates} a \CGrev{Mixed Integer Nonlinear Programming (MINLP)} problem and shows that the level of sparsity can be imposed in the model\revtwoCGblue{. Alternatively,} the variance of the principal components and their sparsity can be jointly maximized in a biobjective framework \cite{CARRIZOSA2014151}.

\subsection{Partial Least Squares}\label{sec:partial-least-squares}
Partial Least Squares (PLS) identifies transformed features $Z_1, \dots, Z_H$ by projecting both the predictors \JSrev{$X$ and their corresponding response $y$ into a new space}, and \revtwoCGblue{this is an approach specific} to regression problems \CGrev{\cite{friedman2001elements}}. \revtwoCG{PLS is particularly viable for problems with a large number of features compared to observations as it aims to identify the latent factors that explain most the variations in the response. PLS corresponds to fitting simple regression models each containing a single predictor variable.} 

The first PLS direction is denoted by $\phi^{1}\in \mathbb{R}^p$ where each component $\phi_{j}^1$ is found by fitting a regression with predictor $X_j$ and response $y$. The first PLS direction points towards the \JSrev{features} that are more strongly related to the response. For computing the second PLS direction, the features vectors $X_1, \dots, X_p$ are first orthogonalized with respect to $Z_1$ (as per the Gram-Schmidt approach), and then individually fitted in simple regression models with response $y$\revtwoCG{. The process is iterated} for all PLS directions $H<p$. The coefficient of the simple regression of $y$ onto each original feature $X_j$ can also be computed as the inner product $\langle y, X_j\rangle$. Similar to PCR, PLS then fits a linear regression model with regressors $Z_1, \dots, Z_H$ and response $y$.

While the principal components directions maximize variance, PLS searches for directions \CGrev{\mbox{$Z_h = \sum_{j=1}^p \phi_{j}^h X_j$}} with both high variance and high correlation with the response. The $h$-th direction $\phi^h$ can be found by solving the optimization problem
\begin{align}
\max_{\phi^h \revtwoCG{\in \mathbb{R}^p}} & \ \CGrev{\text{Corr}(y, X \phi^h)^2\times \text{Var}(X\phi^h)}\\
\mbox{s.t. }& \sum_{j=1}^p (\phi_{j}^h)^2=1 \revtwoCG{,}\\
& {\phi^h}^\top S \phi^l = 0 \quad \forall l=1, \dots, h-1 \revtwoCG{,}\label{eq:PLS-uncorr}
\end{align}
where Corr() indicates the correlation matrix, Var() the \CGrev{variance},
$S$ the sample covariance matrix of $X_j$, and \eqref{eq:PLS-uncorr} ensures that \CGrev{$Z_m$} is uncorrelated with the previous directions \CGrev{$Z_l = \sum_{j=1}^p \phi_{j}^l X_j$}.

\section{Deep Learning}\label{sec:deeplearning}
Deep Learning received a first momentum until the 80s due to universal approximation results \cite{cybenko1989approximation, hornik1991approximation}\revtwoCGblue{. N}eural networks with a single layer with a finite number of units can represent any multivariate continuous function on a compact subset \CGrev{in} $\mathbb{R}^n$ with arbitrary precision. However, the computational complexity required for training Deep Neural Networks (DNNs) hindered their diffusion by late 90s. Starting 2010, the empirical success of DNNs has been widely recognized for several reasons, including the development of advanced processing units, namely GPUs, the advances in the efficiency of training algorithms such as backpropagation, the establishment of proper initialization parameters, and the massive collection of data enabled by new technologies in a variety of domains (e.g., healthcare, supply chain management \cite{TIWARI2018319}, marketing, logistics \cite{WANG201698}, Internet of Things). 

\CGrev{DNNs can be used for the regression and classification tasks discussed in the previous sections, especially when traditional machine learning models fail to capture complex relationships between the input data and the quantitative response, or class, to be learned.}
The aim of this section is to describe the decision optimization problems associated with DNN architectures. \CGrev{To facilitate the presentation, the notation for the common parameters is provided in Table \ref{table:DNN_params}, and an example of fully connected feedforward network is shown in Figure \ref{figure:DNN}.}

\begin{table}[!h]
\begin{tabular}{cp{0.8\linewidth}}
	$\{0, \dots, L\}$ & \textit{layers} indices\revtwoCG{.}\\
	$n^l$ & number of \textit{units}, or \textit{neurons}, in layer $l$\revtwoCG{.}\\
	$\sigma$ & element-wise activation function\revtwoCG{.}\\
	$U(j,l)$ & $j$-th unit of layer $l$\revtwoCG{.}\\
	\CGrev{$W^l \in \mathbb{R}^{n^l \times n^{l+1}}$} & weight matrix for layer \CGrev{$l<L$}\revtwoCG{.} \\
	\CGrev{$b^l \in \mathbb{R}^{n^l}$} & bias vector for layer \CGrev{$l>0$} \revtwoCG{.}\\
	\CGrev{$(X, y)$} & \CGrev{\textit{training} dataset, with observations $x_i$ and responses $y_i, i=1, \dots, n.$}\\
	$x^l$ & \CGrev{output vector of layer $l$ ($l=0$ indicates \textit{input feature} vector, $l>0$ indicates \textit{derived feature} vector)}.	
\end{tabular}\caption{Notation for DNN architectures.}\label{table:DNN_params}
\end{table}


\def\layersep{2.5cm}
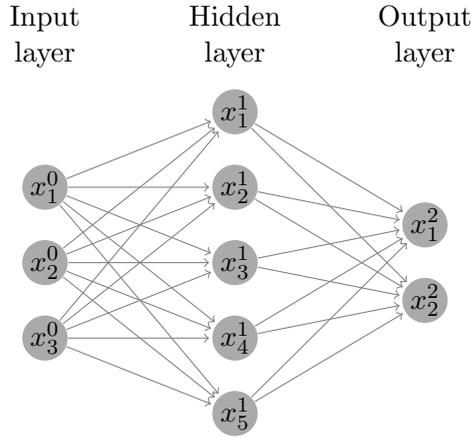
\begin{figure} \centering
\begin{tikzpicture}[shorten >=1pt,->,draw=black!50, node distance=\layersep]
\tikzstyle{every pin edge}=[<-,shorten <=1pt]
\tikzstyle{neuron}=[circle,fill=black!25,minimum size=17pt,inner sep=0pt]
\tikzstyle{input neuron}=[neuron, fill={rgb:black,1;white,2}];
\tikzstyle{output neuron}=[neuron, fill={rgb:black,1;white,2}];
\tikzstyle{hidden neuron}=[neuron, fill={rgb:black,1;white,2}];
\tikzstyle{annot} = [text width=4em, text centered]

\foreach \name / \y in {1,...,3}
\node[input neuron] (I-\name) at (0,-\y -0.5) {$x^0_\y$};

\foreach \name / \y in {1,...,5}
\path[yshift=0.5cm]
node[hidden neuron] (H1-\name) at (\layersep,-\y cm) {$x^1_\name$};

\node[output neuron] (O-1) at (2*\layersep, -2 cm) {$x^2_1$} ; 
\node[output neuron] (O-2)  at (2*\layersep, -3 cm) {$x^2_2$}; 

\foreach \source in {1,...,3}
\foreach \dest in {1,...,5}
\path (I-\source) edge (H1-\dest);

\foreach \source in {1,...,5}
\path (H1-\source) edge (O-1);

\foreach \source in {1,...,5}
\path (H1-\source) edge (O-2);

\node[annot,above of=H1-1, node distance=1cm] (hl) {Hidden layer};
\node[annot,left of=hl] {Input layer};
\node[annot,right of=hl] {Output layer};
\end{tikzpicture}
\caption{\CGrev{Deep Feedforward Neural Network with 3 layers. The input layer has $n^0=3$ units, \revtwoCG{the} hidden layer has $n^1=5$ units and there are $n^2=2$ output units. This is an example of fully connected network, where each neuron in one layer is connected to all neurons in the next layer.
		 Training such network requires to determine weight matrices \revtwoCG{$W^0 \in \mathbb{R}^{3\times5}$, $W^1 \in \mathbb{R}^{5\times2}$}, 	and bias vectors \revtwoCG{$b^1 \in \mathbb{R}^5$}, $b^2 \in \mathbb{R}^{2}.$
}}
	\label{figure:DNN}
\end{figure}

The output vector $x^L$ of a DNN is computed by propagating the information from the input layer to each following layer \CGrev{via the weight matrices $W^l, \ l < L$, the bias vectors $b^l, \ l > 0$, and the activation function $\sigma$}, such that
\begin{align} x^l = \sigma(W^{l-1}x^{l-1} + b^{l-1}) \quad \forall l = 1, \dots, L. \label{eq:computationNN} \end{align}
\CGrev{Activation functions indicate whether a neuron should be activated or not in the network, and are responsible for the capability of DNNs to learn complex relationships between the input and the output.}
\revtwoCGblue{T}he \textit{rectified linear unit} $$ReLU: \mathbb{R}^n \rightarrow \mathbb{R}^n, ReLU(z) = (\max(0,z_1), \dots, \max(0,z_n))$$ is typically one of the preferred options \revtwoCGblue{for activation functions}\CGrev{, mainly because it can be optimized with gradient-based methods \revtwoCGblue{for DNN training}, and tends to produce sparse networks (where not all neurons are activated)}. 

In the context of regression, the components of $x^L$ can directly represent the response values learned.
For a classification problem, the vector $x^L$ corresponds to the \textit{logits} of the classifier. In order to interpret $x^L$ as a vector of class probabilities, functions $F$ such as the logistic sigmoidal or the softmax can be applied \revtwoCGblue{to \mbox{$x^L$ \cite{goodfellow2016deep}}}. The classifier $\mathcal{C}$ modeled by the DNN then classifies an input $x$ with the label correspondent to the maximum activation $\mathcal{C}(x) = \displaystyle\argmax_{i=1,\dots,n^L} F(x_i^L).$

The task of training a DNN consists of determining the weights $W^l$ and the biases $b^l$ that make the model best fit the training data, according to a certain measure of training loss.
\revtwoCG{In multivariate regression with 
\CGrev{$K$} response variables \cite{izenman2008modern, mielke1997permutation}, the training loss $\mathcal{L}$ is typically the sum-of-squared errors $\displaystyle\sum_{k=1}^{K} \displaystyle\sum_{i=1}^{n} (y_{ik} - x^L_k)^2$ where $y_{ik}$ denotes response $k$ corresponding to the $i$-th input vector.}
For classification with \CGrev{$K$} classes, cross-entropy $- \displaystyle\sum_{k=1}^{K} \displaystyle\sum_{i=1}^{n} y_{ik} \log x^L_k$ is preferred. An effective approach to minimize $\mathcal{L}$ is by gradient descent, called \textit{back-propagation} in this setting. Typically, one is not interested in a proven local minimum of $\mathcal{L}$, as this is likely to overfit the training dataset and yield a learning model with a high variance. Similar to the Ridge regression (see Section \ref{sec:regression}), the loss function can include \CGrev{regularization terms, such as }a weight decay term $$\lambda \bigg(\displaystyle\sum_{l=0}^{L-1} \displaystyle\sum_{i=1}^{n_l} (b_{i}^{l})^2 + \displaystyle\sum_{l=0}^{L-1} \displaystyle\sum_{i=1}^{n_l} \displaystyle\sum_{j=1}^{n_{l+1}} (W_{ij}^{l})^2  \bigg),$$
or alternatively a weight elimination penalty term $$\lambda \bigg(\displaystyle\sum_{l=0}^{L-1} \displaystyle\sum_{i=1}^{n_l} \frac{(b_{i}^{l})^2}{1+(b_{i}^{l})^2} + \displaystyle\sum_{l=0}^{L-1} \displaystyle\sum_{i=1}^{n_l} \displaystyle\sum_{j=1}^{n_{l+1}} \frac{(W_{ij}^{l})^2}{1+(W_{ij}^{l})^2}  \bigg).$$
\CGrev{Weight decay limits the growth of the weights, which \revtwoCGblue{speeds up the training via backpropagation}, and has been shown to limit overfitting (\revtwoCG{see} \cite{poggio2017theory} for a discussion about overfitting in Neural Networks).}

The aim of this section is to present the optimization models that are used in DNN \CGrev{for feedforward architectures. Several other neural network architectures have been investigated in deep learning \cite{goodfellow2016deep}. In particular, Convolutional Neural Networks (CNN) \cite{lecun1989generalization} have been successfully adopted for processing data with a grid-like topology, such as images \cite{krizhevsky2012imagenet}, videos \cite{karpathy2014large}, and traffic analytics \cite{wang2018enhancing}. In CNN, the output of layers is obtained via convolutions (instead of the matrix multiplication in feedforward networks), and pooling operations on nearby units (such as average or maximum operators).}
\CGrev{In the remainder of the section, mixed integer} programming models for DNN training are introduced in Section~\ref{sec:dnn_training_mip}\CGrev{, and ensemble approaches with multiple activation functions are discussed in Section~\ref{sec:ensemble_approaches}}. 

\subsection{\CGrev{Mixed Integer} Programming for DNN Architectures} \label{sec:dnn_training_mip}

Motivated by the considerable improvements of  mixed integer programming solvers, a natural question is how to model a \revtwoCG{trained} DNN as a MIP. In \cite{fischetti2018deep}, DNNs with \textit{ReLU} activation \begin{align} x^l = ReLU(W^{l-1}x^{l-1} + b^{l-1}) \quad \forall l = 1, \dots, L \label{eq:computationNNReLU} \end{align}
are modeled \CGrev{as a MIP with decision variables $x^l$ expressing the output vector of layer $l$, $l>0$ and $l^0$ is the input vector. To express \eqref{eq:computationNNReLU} explicitly, each unit $U(j,l)$ of the DNN is associated with binary activation variables $z_j^l$, and continuous slack variables $s_j^l$.} The following mixed integer linear problem is proposed
\begin{align}
\min \ & \sum_{l=0}^L\sum_{j=1}^{n_l} c_j^l x_j^l + \sum_{l=1}^L\sum_{j=1}^{n_l} \gamma_j^l z_j^l \label{DNN_Fischetti_obj}\\
\mbox{s.t. } & \sum_{i=1}^{n_{l-1}} w_{ij}^{l-1} x_i^{l-1} + b_j^{l-1} = x_j^l - s_j^l \qquad \forall l = 1, \dots, L, j = 1, \dots, n_l\revtwoCG{,} \label{DNN_Fischetti_prop}\\
& \CGrev{x_j^l \leq (1-z_j^l)M_x^{j, l}} \qquad \forall l = 1, \dots, L, j = 1, \dots, n_l\revtwoCG{,} \label{act_z_1}\\
& \CGrev{s_j^l \geq z_j^l M_s^{j,l}}  \qquad \forall l = 1, \dots, L, j = 1, \dots, n_l\revtwoCG{,} \label{act_z_2}\\
& \CGrev{0} \leq x_j^l \leq ub^l_j \qquad \forall l = 1, \dots, L, j = 1, \dots, n_l \revtwoCG{,}
\label{cons:DNN-xk}\\
& \CGrev{0} \leq s_j^l \leq \overline{ub}^l_j \qquad \forall l = 1, \dots, L, j = 1, \dots, n_{\revtwoCG{l}} \revtwoCG{,} \label{DNN_Fischetti_fin}
\end{align}
where \CGrev{$M_x^{j, l}, M_s^{j, l}$ are suitably large constants}. \revtwoCG{We note that since the DNN is trained, the weights $w_{ij}^{l}$ and bias $b_j^{l}$ are fixed parameters}. Depending on the application, different activation weights $c_j^l$ and activation costs $\gamma_j^l$ can also be used for each $U(j,l)$. \CGrev{If known, upper bound $ub_j^l$ can be enforced on the output $x_j^l$ of unit $U(j,\revtwoCG{l})$ via constraints \eqref{cons:DNN-xk}, and slack $s_j^l$ can be bounded by \revtwoCG{$\overline{ub}^l_j$} via constraints \eqref{DNN_Fischetti_fin}.
}



The proposed MIP is feasible for every input vector $x^0$ since it computes the activation in the subsequent layers. \CGrev{Constraints \eqref{act_z_1} and \eqref{act_z_2} are known to have a weak continuous relaxation, and the tightness of the chosen constants (bounds) is crucial for their effectiveness. Several optimization solvers can directly handle such kind of constraints as indicator constraints~\cite{Bonami2015}.}
In \cite{fischetti2018deep}, a bound-tightening strategy to reduce the computational times is proposed and the largest DNN tested with this approach is a $5$-layer DNN with $20+20+10+10+10$ internal units.

\CGrev{Problem \eqref{DNN_Fischetti_obj}--\eqref{DNN_Fischetti_fin} can model several tasks in Deep Learning, other than the computation of quantitative responses in regression, and of classification. Such tasks include}
\begin{itemize}
	\item \textit{Pooling operations:} The average and the maximum  operators
	\begin{align*}
	& \textit{Avg}(x^l) = \frac{1}{n^l} \sum_{i=1}^{n^l} x_i^l \revtwoCG{,} \\
	& \textit{Max}(x^l) = \max(x_1^l, \dots, x_{n^l}^l),
	\end{align*}
	\CGrev{can be incorporated in the hidden layers}. In the case of max pooling operations, additional indicator constraints are required. \revtwoCGblue{For example, a}\CGrev{verage and maximum operators are often used in CNNs, as mentioned earlier in \Cref{sec:deeplearning}.}
	\item \textit{Maximizing the unit activation:} By maximizing the objective function \eqref{DNN_Fischetti_obj}, one can find input examples $x^0$ that maximize the activation of the units. This may be of interest in applications such as the visualization of image features.
	\item \textit{Building crafted adversarial examples:} Given an input vector $x^0$ labeled as $\chi$ by the DNN, the search for perturbations of $x^0$ that are classified as $\chi' \neq \chi$ (\textit{adversarial examples}), can be conducted by adding conditions on the activation of the final layer $L$ and minimizing the perturbation. \cg{In \cite{fischetti2018deep}, such conditions are actually restricting the search for adversarial examples and the resulting formulation does not guarantee an adversarial solution nor can prove that no adversarial examples exist.} Adversarial learning \revtwoCGblue{is the objective of the discussion in} \Cref{sec:adv_learning}. 
	\item \textit{Training:} In this case, the weights and biases are decision variables. The resulting bilinear terms in \eqref{DNN_Fischetti_prop} and the considerable number of decision variables in the formulation limit the applicability of \eqref{DNN_Fischetti_obj}--\eqref{DNN_Fischetti_fin} for DNN training. 
\end{itemize}
	
 Another attempt in modelling DNNs via MIPs is provided by \cite{khalil2018combinatorial}, in the context of Binarized Neural Networks (BNNs). BNNs are characterized by having binary weights $\{-1, +1\}$ and by using the sign function for neuron activation  \cite{courbariaux2015binaryconnect}. 
 In \cite{khalil2018combinatorial}, a MIP is proposed for finding adversarial examples in BNNs by maximizing the difference between the activation of the targeted label $\chi'$ and the predicted label $\chi$ of the input $x^0$, in the final layer (namely, $\max x_{\chi'}^L - x_\chi^L$).
  \comments{Contrary to \cite{fischetti2018deep}, the MIP of \cite{khalil2018combinatorial} does not impose limitations on the search of adversarial examples, apart from the perturbation quantity. In terms of optimality criterion however, searching for the proven largest misclassified example is different from finding a targeted adversarial example. Furthermore, while there is interest in minimally perturbed adversarial examples, suboptimal solutions corresponding to adversarial examples (i.e., $x_{\chi'}^L \geq x_\chi^L$) may have a perturbation smaller than that of the optimal solution.} \CGrev{Recently, \cite{icartetraining} investigated a hybrid constraint programming/mixed integer programming method to train BNNs. Such model-based approach provides solutions that generalize better than those found by the largely adopted training solvers, such as gradient descent, especially for small datasets. \revtwoJS{We note that methods such as gradient descent can usually only guarantee local optimality (unless early stopping takes place).}}
  

Besides \cite{fischetti2018deep}, other MIP frameworks have been proposed to model certain properties of neural networks in a bounded input domain. 
In \cite{cheng2017maximum}, the problem of computing maximum perturbation bounds for DNNs  is formulated as a MIP, where indicator constraints and disjunctive constraints are modeled using constraints with big-M coefficients \cite{grossmann2002review}. \comments{The maximum perturbation bound is a threshold such that the perturbed input may be classified correctly with a high probability. A restrictive misclassification condition is added when formulating the MIP.  Hence, the infeasibility of the MIP does not certify the absence of adversarial examples.} In addition to the $ReLU$ activation, the $tan^{-1}$ function is also considered by introducing quadratic constraints and several heuristics are proposed to solve the resulting problem.
In \cite{tjeng2017verifying}, a model to formally measure the vulnerability to adversarial examples is proposed (the concept of vulnerability of neural networks is discussed in more details in Sections \ref{sec:targeted-attacks} and \ref{sec:untargeted-attacks}).
 A tight formulation for the resulting nonlinearities and a novel presolve technique are introduced to limit the number of binary variables and improve the numerical conditioning. \comments{However, the misclassification condition of adversarial examples is not explicitly defined but is rather left in the form \enquote{different from} and not explicitly modeled using equality/inequality constraints.}
 In \cite{serra2017bounding}, the aim is to count or bound the number of linear regions that a piecewise linear classifier represented by a DNN can attain. Assuming that the input space is bounded and polyhedral, the DNN is modeled as a MIP. \comments{The contributions of adopting a MIP framework in this context are limited, especially in comparison with the computational results achieved in \cite{montufar2017notes}.} 

MIP frameworks can also be used to formulate the verification problem for neural networks as a satisfiability problem. In \cite{katz2017ReLUplex}, a satisfiability modulo theory solver is proposed based on an extension of the simplex method to accommodate the $ReLU$ activation functions. In \cite{BunTTKK18}, a branch-and-bound framework for verifying piecewise-linear neural networks is introduced. For a recent survey on the approaches for automated verification of NNs, the reader is referred to \cite{leof2018automver}. 

\subsection{Activation Ensembles}\label{sec:ensemble_approaches}
Another research direction in neural network architectures investigates the possibility of adopting multiple activation functions inside the layers of a neural network, to increase the accuracy of the classifier. 
Some examples in this framework are given by the \textit{maxout} units \cite{goodfellow2013maxout}, returning the maximum of multiple linear affine functions, and the \textit{network-in-network} paradigm \cite{lin2013network} where the $ReLU$ activation function is replaced by fully connected network. 
In \cite{agos2014}, adaptive piecewise linear activation functions are learned when training \CGrev{each neuron}. Specifically, for each unit $i$ and value $z$, activation \CGrev{$\sigma_i(z)$} is considered as
\begin{align}
\sigma_i(z) = \max(0, z) + \sum_{s=1}^S a_i^s \max (0, -z + b_i^s), \label{eg:agos}
\end{align}
where the number of hinges $S$ is a hyperparameter to be fixed \revtwoCGblue{before training}, while the variables $a_i^s, b_i^s$ have to be learned. \CGrev{Functions \revtwoCG{$\sigma_i$} generalize the   
$ReLU$ function (first term of \eqref{eg:agos}), and can approximate a class of continuous piecewise-linear functions, for large enough $S$ \cite{agos2014}.} 

In a more general perspective, ensemble layers are proposed in \cite{2017arXiv170207790H} to consider multiple activation functions in a neural network. \CGrev{The idea is to embed a family of activation functions $\{\Phi^1, \dots, \Phi^m\}$ and let the network itself choose the magnitude of their activation \CGrev{for each neuron $i$} during the training.} To promote relatively equal contribution to learning, the activation functions need to be scaled to the interval $[0,1]$. \revtwoCGblue{In order to} measure the impact of the activation in the neural network, each function \CGrev{$\Phi^j$} is associated with a continuous variable $\alpha^j$. 
\CGrev{The resulting activation $\sigma_i$ for neuron $i$ is then given by}
\begin{align}
& \CGrev{\sigma_i(z) = \sum_{j=1}^{m}  \alpha_i^j \cdot \frac{\Phi^j(z) - \displaystyle\min_{x \in X}(\Phi^j(z_{x, i}))}
	{\displaystyle\max_{x \in X}(\Phi^j(z_{x, i})) - \displaystyle\min_{x \in X}(\Phi^j(z_{x, i})) + \epsilon} }\revtwoCG{,} \label{act-ens}
\end{align}
where \CGrev{$z_{x, i}$} is the \CGrev{output of neuron $i$} associated with training example \CGrev{$x$}, \CGrev{$X$ is the set of training observations}, and $\epsilon$ is a small tolerance. 
Equation \eqref{act-ens} \CGrev{is a weighted sum of the scaled $\Phi^j$ functions, which is integrated in the training of the DNN architecture. The $\min$ and $\max$ in \eqref{act-ens} can be approximated on a minibatch of observations in $X$, in the testing phase.}   
 \CGrev{In order to impose the selection of functions $\Phi^j$, the} magnitude of the weights $\alpha^j$ is limited in a projection subproblem, where for each neuron the network should choose an activation function and \CGrev{therefore all $\alpha^j$ should sum to $1$}. If $\hat{\alpha_j}$ are the weight values obtained by gradient descent while training, then the projected weights are found by solving \CGrev{the convex quadratic programming problem}
\begin{align}
\min \ & \sum_{j=1}^{m} \frac{1}{2} (\alpha^j - \hat{\alpha}^j)^2 \label{adv_ens_proj:obj}\\
\mbox{s.t. } & \sum_{j=1}^{m} \alpha^j = 1\revtwoCG{,} \label{adv_ens_proj:alpha_lim}\\
& \CGrev{\alpha^j} \geq 0 \; \forall j=1, \dots, m, \label{adv_ens_proj:alpha_bound}
\end{align}
\CGrev{which can be solved in closed form via the Karush-Kuhn-Tucker (KKT) conditions.}

\section{\CGrev{Adversarial Learning}}\label{sec:adv_learning}

\CGrev{Despite the wide adoption of Machine Learning models in real-world applications, their integration into safety and security related use cases still necessitates thorough evaluation and research. A large number of contributions in the literature pointed out the dangers caused by perturbed examples, also called \textit{adversarial examples}, causing classification errors \cite{biggio2010multiple, szegedy2013intriguing}.} Malicious attackers can thus exploit security falls in a general classifier. In case the attacker has a perfect knowledge of the classifier's architecture (i.e., the result of the training phase), then a \textit{white-box} attack can be performed. \textit{Black-box} attacks are instead performed without full information of the classifier. The \CGrev{interest in} adversarial examples is also motivated by the transferability of the attacks to different trained models \cite{kurakin2016adversarial, tramer2017ensemble}. 
Adversarial learning then emerges as a framework to devise vulnerability attacks for classification models \cite{lowd2005adversarial}.

From a mathematical perspective, such security issues have been formerly expressed via min-max approaches where the learner's and the attacker's loss functions are antagonistic \cite{dekel2010learning, globerson2006nightmare, lanckriet2002robust}. Non-antagonistic losses are formulated as a Stackelberg equilibrium problem involving a bilevel optimization formulation \cite{bruckner2011stackelberg}, or in a Nash equilibrium approach \cite{bruckner2012static}. These theoretical frameworks rely on the assumption of expressing the actual problem constraints in a game-theory setting, which is often not a viable option for real-life applications.

\CGrev{The search for adversarial examples can also be used to evaluate the efficiency of Generative Adversarial Networks (GANs) \cite{goodfellow2014generative}. A GAN is a minmax two-player game where a generative model $G$ tries to reproduce the training data distribution and a discriminative model $D$ estimates the probability of detecting samples coming from the true training distribution, rather than $G$. The game terminates at a saddle point\revtwoCGblue{,} \revtwoCG{which} is a minimum with respect to a player's strategy and a maximum for the other player's strategy. Discriminative networks can be affected by the presence of adversarial examples because the specific inputs to the classification networks are not considered in GANs training.}

\revtwoCGblue{A}dversarial attacks on the test set can be conducted in a targeted or untargeted fashion \cite{carlini2017towards}. In the targeted setup, the attacker aims to achieve a classification with a chosen target class \revtwoCGblue{(discussed in \Cref{sec:targeted-attacks})}, while the untargeted misclassification is not constrained to \CGrev{achieve a} specific class \revtwoCGblue{(\Cref{sec:untargeted-attacks}).} The robustness of DNNs to adversarial attacks is discussed in \Cref{sec:robust}\CGrev{. Finally, data poisoning attacks are described in \Cref{sec:data-poisoning}. While the majority of the cited papers of the present section refer to DNN applications, adversarial learning can, in general, be formulated for classifiers \revtwoCGblue{with quantitative classes}, such as those discussed in \Cref{sec:classification}.}

\subsection{Targeted attacks}\label{sec:targeted-attacks}
\CGrev{
\revtwoCG{Given a neural network classifier $f: \psi \subset \mathbb{R}^{p} \rightarrow \Upsilon$, an input $x \in \psi$ with label $y \in \Upsilon$, and a target label $y' \in \Upsilon$, a targeted attack consists of a perturbation $r$ such that $f(x+r)=y'$.}
}
This corresponds to finding an input \enquote{close} to $x$, which is misclassified by $f$. Clearly, if the target $y'$ coincides with $y$, the problem has the trivial solution $r=0$ and no misclassification takes place.

The minimum adversarial problem for targeted attacks consists of finding a perturbation $r$ by solving
\begin{align}
\min_{r \in \mathbb{R}^{p}} \ & \lVert r \rVert_2 \label{min_dist:obj}\\
\mbox{s.t. } & f(x+r) = y' \revtwoCG{,} \label{min_dist:target}\\
& x+r \in \revtwoCG{\psi}. \label{min_dist:valid}
\end{align}
The condition \eqref{min_dist:valid} ensures that the perturbed example $x+r$ belongs to the set of admissible inputs. 
The difficulty of solving problem \eqref{min_dist:obj}--\eqref{min_dist:valid} to optimality depends on the complexity of the classifier $f$\revtwoCG{, and the set $\psi$ of feasible inputs. In general, it is computationally challenging to find an optimal solution to the problem, especially in the case of neural networks.}

\revtwoCG{For classification of normalized images with binary pixel values, \cite{szegedy2013intriguing} introduces the box-constrained approximation}
 \CGrev{
\begin{align}
\min_{r \in \mathbb{R}^{p}} \ & c |r| + \mathcal{L} (x+r, y') \label{min_dist_L-BFGS:obj}\\
\revtwoCG{\mbox{s.t. }} & x+r \in \revtwoCG{[0,1]^{p}}, \label{min_dist_L-BFGS:valid}
\end{align}
\revtwoCG{where $\mathcal{L} :
	\psi \times \Upsilon \rightarrow \mathbb{R}^+$ denotes the loss function for training $f$ (e.g., cross-entropy).}
The approximation is exact for convex loss functions, and can be solved via a line search algorithm on $c>0$.
For a fixed $c$, the formulation can be tackled by the box-constrained version of the Limited-memory Broyden–Fletcher–Goldfarb–Shanno (L-BFGS) method \cite{byrd1995limited}.}
In \cite{gu2014towards}, $c$ is fixed such that the perturbation is 
minimized on a sufficiently large subset \CGrev{$X'$} of data points, and the mean prediction error rate of $f(x_i + r_i),  x_i \in X'$ is greater than a threshold. 
In \cite{carlini2017towards}, the \revtwoCG{$2$-norm in \eqref{min_dist:obj}} is generalized to the \JSrev{$l$-norm} with $l\in \{0, 2, \infty\}$ and an alternative formulation is introduced which includes functions~$\mathcal{F}$ in the objective where $f(x+r)=y'$ is satisfied if and only if $\mathcal{F}(x+r) \leq 0$. The equivalent formulation is then
\begin{align}
\min_{r \in \mathbb{R}^{p}} \ & \lVert r \rVert_l + \Lambda \mathcal{F}(x+r) \label{min_dist_nobox:obj}\\
\revtwoCG{\mbox{s.t. }} & x+r \in \revtwoJS{\psi}, \label{min_dist_nobox:valid}
\end{align}
where $\Lambda$ is a constant that can be determined by binary search such that the solution $r^*$ satisfies the condition $\mathcal{F}(x+r^*) \leq 0$.
For the case where \eqref{min_dist_nobox:valid} are box constraints similar to \eqref{min_dist_L-BFGS:valid}, the authors propose strategies for applying optimization algorithms such as Adam \cite{Adam}. Novel classes of attacks are identified for the considered metrics.

\subsection{Untargeted attacks}\label{sec:untargeted-attacks}

\CGrev{
In untargeted attacks, one searches for adversarial examples $x'$ close to the original input $x$ with label $y$ for which the classified label $y'$ of $x'$ is different from $y$, without targeting a specific label for $x'$.
} Given that the only aim is misclassification,
untargeted attacks are deemed less powerful than the targeted counterpart, and received less attention in the literature.

\CGrev{
A mathematical formulation for finding minimum adversarial distortion for untargeted attacks is proposed in \cite{tjeng2017verifying}. Assuming that the output values of classifier $f$ are expressed by the functions $f_{y'}$ associated with labels \revtwoCG{$y' \in \Upsilon$} (i.e., $f_{y'}$ are the scoring functions), and a distance metric $d$ is given, then a \revtwoCGblue{minimum} perturbation $r$ for an untargeted attack is found by solving
}
\begin{align}
\min_{r\in \mathbb{R}^{p}} \ & d(r) \label{adv_dist:obj}\\
\mbox{s.t. } & \argmax_{y' \in \Upsilon} \{f_{y'}(x+r)\} \neq y \revtwoCG{,} \label{adv_dist:argmax}\\
& x+r \in \psi. \label{adv_dist:valid}
\end{align}
This formulation can easily accommodate targeted attacks in a set \CGrev{$T \not\owns y$} by replacing \eqref{adv_dist:argmax} with $\argmax_{y'} \{f_{y'}(x+r)\} \in T$. The most commonly adopted metrics in the literature are the \revtwoCG{$1$, $2$, and $\infty$-norm, which can all be expressed with continuous variables, as shown in \cite{tjeng2017verifying}}. The \CGrev{$2$-norm} makes the objective function of the outer-level optimization problem quadratic. 

	\CGrev{In order to express the logical constraint \eqref{adv_dist:argmax} in a mathematical programming formulation, }
	problem \eqref{adv_dist:obj}--\eqref{adv_dist:valid} can be cast as the bilevel optimization problem 
	\begin{align}
	\min_{r\in \mathbb{R}^{p}, z\in \Upsilon} \ & d(r) \label{adv_dist_bi:obj}\\
	\mbox{s.t. } & z - y \leq - \epsilon + M s \revtwoCG{,}\label{adv_dist_bi:neq1}\\
	& z - y \geq \epsilon -(1-s) M \revtwoCG{,}\label{adv_dist_bi:neq2}\\
	& z \in \argmax_{y' \in \Upsilon} \{f_{y'}(x+r)\} \revtwoCG{,}\label{adv_dist_bi:argmax} \\
	& x+r \in \psi,\label{adv_dist_bi:fin}\\
	& s \in \{0,1\},
	\end{align}
	\CGrev{where $\epsilon>0$ is a small constant, $z$ is a decision variable representing the classified label, $M$ is a big-M coefficient, and $s$ is a binary variable that enforces one of the constraints }
	\eqref{adv_dist_bi:neq1}--\eqref{adv_dist_bi:neq2} which express the condition of misclassification $z \neq y$. 
\CGrev{The complexity of the inner-level optimization problem is dependent on the scoring functions. 
Given that the upper-level feasibility set $\psi$ is typically continuous and the lower-level variable $y'$ ranges on a discrete set, the problem is in fact a continuous discrete bilevel programming problem \cite{fanghanel2009bilevel} with convex quadratic function \cite{edmunds1992algorithm}, which requires dedicated reformulations or approximations \revtwoCG{\cite{chen1995nonlinear, gumucs2001global, jan1994nonlinear}}. 
}
	
	
	We introduce an alternative mathematical formulation for finding untargeted adversarial examples \CGrev{satisfying condition \CGrev{\eqref{adv_dist:argmax}}}. \revtwoCG{A perturbed input $x' = x +r$ for a sample $x$ classified with label $y \in \Upsilon$ is an untargeted adversarial example if the classified label of $x'$ is different from $y$. This condition is equivalent to
	\begin{align}
	& \exists\ y' \in \Upsilon\setminus\{y\}\ \text{s.t.} \ f_{y'}(x') > f_{y}(x'). \label{untargeted}
	\end{align}}
	Condition \eqref{untargeted} is an existence condition, which can be formalized by \CGrev{introducing} the functions $\tilde{\sigma}_{y'}(r) = \CGrev{ReLU}(f_{y'}(x+r) -  f_{y}(\revtwoCG{x+r}))$, $y' \in \Upsilon\setminus\{y\}$, and the condition
	\begin{align}
	& \sum_{y' \in \Upsilon\setminus\{y\}} \tilde{\sigma}_{y'}(r) > \CGrev{\nu} \label{untargeted_ReLU},
	\end{align}
	where \CGrev{parameter $\nu > 0$ enforces that at least one $\tilde{\sigma}_{y'}$ function has to be activated for a perturbation $r$. Therefore, untargeted adversarial examples can be found from formulation \eqref{adv_dist:obj}--\eqref{adv_dist:valid} by replacing condition \eqref{adv_dist:argmax} with the linear condition \eqref{untargeted_ReLU} and	
	adding \revtwoCGblue{$\lvert \Upsilon \rvert-1$} functions $\tilde{\sigma}_{y'}(r)$}. The complexity of this approach depends on the scoring functions $f_{y'}$. The extra \CGrev{$ReLU$} functions \CGrev{$\tilde{\sigma}$} can be expressed as a \CGrev{mixed integer formulation} as done in problem \eqref{DNN_Fischetti_obj}--\eqref{DNN_Fischetti_fin}.



\subsection{Adversarial robustness}\label{sec:robust}

\CGrev{Another interesting line of research motivated by adversarial learning deals with adversarial training, which consists of techniques to make a neural network robust to adversarial attacks. The problem of} measuring \revtwoCG{the} robustness of a neural network is formalized in \cite{bastani2016measuring}. The \textit{pointwise robustness} evaluates if the classifier $f$ on $x$ is robust for \enquote{small} perturbations.
Formally, $f$ is said to be $(x, \epsilon)$-robust if
\begin{align}
y' = y, \; \forall x' \, \text{s.t.} \, \lVert x'-x \rVert_\infty \leq \epsilon. \label{robustness}
\end{align}
Then, the pointwise robustness $\rho(f, x)$ is the minimum $\epsilon$ for which $f$ fails to be $(x, \epsilon)$-robust:
\begin{align}
\rho(f, x) = \text{inf}\{ \epsilon \geq 0 \; | \; f \text{ is not }  (x, \epsilon)\text{-robust} \}. \label{pointwise}
\end{align}
As detailed in \cite{bastani2016measuring}, $\rho$ is computed by expressing \eqref{pointwise} as a constraint satisfiability problem. \CGrev{By imposing a bound on the perturbation, an estimation of the pointwise robustness can be performed by solving a MIP \cite{cheng2017maximum}.
}



\CGrev{A widely known defense technique is to augment the training data with adversarial examples; this however does not offer robustness guarantees on novel kinds of attacks.
The adversarial training of neural network via robust optimization is investigated in \cite{madry2017towards}.} 
In this setting, the goal is to train a neural network to be resistant to all attacks belonging to a certain class of perturbations. Particularly, the adversarial robustness with a saddle point (min-max) formulation is studied in \cite{madry2017towards} which is obtained by augmenting the Empirical Risk Minimization paradigm. 

\CGrev{Let $\theta \in \mathbb{R}^p$ be the set of model parameters to be learned, and $\mathcal{L}(\theta; x, y)$ be the loss function considered in the training phase (e.g., the cross-entropy loss) \CGrev{for training examples $x\in X$ and labels $y\in\Upsilon$}, and let $\mathcal{S}$ be the set of allowed perturbations (e.g., an $L_{\infty}$ ball).} 
The aim is to minimize the worst expected adversarial loss on the set of inputs perturbed by $\mathcal{S}$
\begin{align}
\min_{\theta}  \mathbb{E}_{(x, y)} \Bigl[ \max_{r \in \mathcal{S}} \mathcal{L}(\theta; x+r, y) \Bigr],     \label{saddle:obj}
\end{align}
\CGrev{where the expectation value is computed on the distribution of the training samples.}
The saddle point problem \eqref{saddle:obj} is viewed as the composition of an \textit{inner maximization} and an \textit{outer minimization} problem. The inner problem corresponds to attacking a trained neural network by means of the perturbations $\mathcal{S}$. The outer problem deals with the training of the classifier in a robust manner. The importance of formulation \eqref{saddle:obj} stems both from the formalization of adversarial training and \CGrev{from the quantification of the robustness given by the objective function value on the chosen class of perturbations}. To find solutions to \eqref{saddle:obj} in a reasonable time, the structure of the local minima of the loss function can be explored.

\CGrev{
Another robust training approach consists of optimizing the
model parameters $\theta$ with respect to worst-case data \cite{shaham2018understanding}. This is formalized by introducing a perturbation set $\mathcal{S}_x$ for each training example $x$. The aim is then to optimize
\begin{align}
\min_{\theta}  \sum_{x\in X} \max_{r \in \mathcal{S}_x} \mathcal{L}(\theta; x+r, y).     \label{shaham:obj}
\end{align}
An alternating ascent and descent steps procedure can be used to solve \eqref{shaham:obj} with the loss function approximated by the first-order Taylor expansion around the training points. 
}

\subsection{Data Poisoning}\label{sec:data-poisoning}
\CGrev{
A popular class of attacks for decreasing the training accuracy of classifiers is that of data poisoning, which was first studied for SVMs \cite{biggio2012poisoning}. A data poisoning attack consists of hiding corrupted, altered or noisy data in the training dataset.}
In \cite{SteinhardtKL17}, 
worst-case bounds on the efficacy of a class of causative data poisoning attacks are studied. The causative attacks \cite{barreno2010security} proceed as follow:
\begin{itemize}
	\item \CGrev{a clean training dataset $\Gamma_C$ with $n$ data points drawn by a data-generating distribution is generated}
	\item the attacker adds malicious examples $\Gamma_M$ \CGrev{to $\Gamma_C$, to let the defender (learner) learn a bad model}
	\item the defender learns model \CGrev{with parameters} $\hat{\theta}$ from the full dataset $\Gamma = \Gamma_C \cup \Gamma_M$, reporting a test loss $\mathcal{L} (\hat{\theta})$.
\end{itemize}

\CGrev{Data poisoning can be viewed as a game between the attacker and the defender players, where the defender wants to minimize $\mathcal{L} (\hat{\theta})$, and the attacker \revtwoCG{seeks} to maximize it.} As discussed in \cite{SteinhardtKL17}, data sanitization defenses to limit the increase of test loss $\mathcal{L} (\hat{\theta})$ include two steps: (i) data cleaning (e.g., removing outliers which are likely to be poisoned examples), to produce a feasible dataset 
\revtwoCGblue{$\Gamma'$}, and (ii) minimizing a margin-based loss on the cleaned dataset $\Gamma \cap \Gamma'$. The learned model is then $\hat{\theta} = \argmin_{\theta \in \Theta} \mathcal L(\theta; \Gamma \cap \Gamma')$.
Poisoning attacks can also be performed in \textit{semi-online} or \textit{online} fashion, where training data is processed in a streaming manner, and not in fixed batches \CGrev{(i.e., offline)}. In the semi-online context, the attacker can modify part of the training data stream so as to maximize the classification loss, and the evaluation of the objective (loss) is done only at the end of the training. In the fully-online scenario, the classifier is instead updated and evaluated during the training process. 
In~\cite{WangPoisoning18}, a white-box attacker's behavior in online learning for a linear classifier $w^T x$ \CGrev{(e.g., SVM} with binary labels $y \in \{-1, +1\}$\revtwoCGblue{)} is formulated. \CGrev{The attacker knows the order in which the training data is processed by the learner.} The data stream $S$ arrives in $T$ instants ($S=\{S_1, \dots, S_T\}$\CGrev{, with $S_t = (X_t, y_t)$}) and the classification weights are updated using an online gradient descent algorithm~\cite{zinkevich2003online} such that 
$w_{t+1} = w_t - \eta_t (\nabla \mathcal{L}(w_t, (x_t, y_t))) + \nabla \Omega(w_t),$
where $\Omega$ is a regularization function, $\eta_t$ is the step length of the iterate update, and $\mathcal{L}$ is a convex loss function. \CGrev{Let \revtwoCGblue{$\Gamma^{,}_T$} be the cleaned dataset at time $T$ \CGrev{(which can be obtained, for instance, via the sphere and slab defenses)},  $U$ be a given upper bound on the number of changed examples in $\Gamma$ due to data sanitization, $g$ be the attacker's objective (e.g., classification error on the test set), $|\cdot|$ be the cardinality of a set.
}
The semi-online attacker \CGrev{optimization problem} can then be formulated as
\begin{align}
\max_{S \in \Gamma'_T} & \ g(w_T) & \label{wang_poisoning:obj}\\
\mbox{s.t. } &  | \{S \setminus \Gamma \}| \leq U, &  \label{wang_poisoning:attacker}\\
& w_t = w_0 - \sum_{\tau=0}^{t-1} \eta_{\tau} (\nabla \mathcal{L}(\omega_{\tau}, S_\tau) + \nabla \mathcal{L}(w_\tau)), 1 \leq t \leq T. & \label{wang_poisoning:updated}
\end{align}

Compared to the offline case, \CGrev{the weights $w_t$ to be learned are} a complex function of the data stream $S$, which makes the gradient computation more challenging and the KKT conditions do not hold. The optimization problem can be simplified by considering a convex surrogate for the objective function, given by the logistic loss. In addition, the expectation is conducted over a separate validation dataset and a label inversion procedure is implemented to cope with the multiple local maxima of the classifier function. The fully-online case can also be addressed by replacing objective \CGrev{\eqref{wang_poisoning:obj}} with $\displaystyle\sum_{t=1}^{\revtwoCGblue{T}} g(w_t).$

\section{Emerging Paradigms}\label{sec:emerging}

\subsection{Machine Teaching}
\CGrev{In all Machine Learning tasks discussed so far, the size of the training set of the machine learning models has been considered as a hyperparameter.} The Teaching Dimension problem identifies the minimum size of a training set to correctly teach a model \cite{GOLDMAN199520, Shinohara1991}. The teaching dimension of linear learners, such as Ridge regression, SVM, and logistic regression has been recently discussed in \cite{liu2016teaching}. 
With the intent to generalize the teaching dimension problem to a variety of teaching tasks, \cite{zhu2015machine} and \cite{zhu2018overview} provide the \textit{Machine Teaching} framework. Machine Teaching is essentially an inverse problem to Machine Learning. While in a learning task, the \CGrev{training dataset $\Gamma = (X, y)$} is given and the model parameters $\theta=\theta^*$ have to be determined, the role of a teacher is to let a learner approximately learn a given model $\theta^*$ by providing a proper set $\Gamma$ of training examples (also called \textit{teaching dataset} in this context). A Machine Teaching task requires to select: 
\begin{enumerate*}
    \item[i)] a Teaching Risk TR expressing the error of the learner, with respect to model $\theta^*$; 
    \item[ii)] a Teaching Cost TC expressing the convenience of the teaching dataset, from the prospective of the teacher, weighted by a \CGrev{regularization factor $\lambda$}; 
    \item[iii)] a learner L.
\end{enumerate*}

Formally, machine teaching can be cast as a bilevel optimization problem
\begin{align}
\min_{\Gamma, \theta} \ & \text{TR}(\theta) + \lambda \text{TC}(\Gamma)  \label{general_machine_teaching:upper}\\
\mbox{s.t. } & \theta = \text{L}(\Gamma). \label{general_machine_teaching:lower}
\end{align}
The upper optimization is the teacher's problem and the lower optimization $L(\Gamma)$ is the learner's machine learning problem. 
The teacher is aware of the learner, which could be a classifier (such as those of \Cref{sec:classification}) or a deep neural network. Machine teaching encompasses a wide variety of applications, such as data poisoning attacks
, computer tutoring systems, and adversarial training.

\CGrev{Problem \eqref{general_machine_teaching:upper}--\eqref{general_machine_teaching:lower} is, in general, challenging to \revtwoCG{solve.  However,} for certain convex learners, one can replace the lower problem by the corresponding KKT conditions, and reduce the problem to a single level formulation.} The teacher is typically optimizing over a discrete space of teaching sets, hence, for some problem instances, the submodularity properties of the problem \revtwoCGblue{may be explored}. For problems with a small teaching set, it is possible to formulate the teaching problem as a \CGrev{mixed integer} nonlinear program. The computation of the optimal training set remains, in general, an open problem, and is especially challenging in the case where the learning algorithm does not have a closed-form solution with respect to the training set \cite{zhu2015machine}. 

\CGrev{The minimization of teaching cost can be directly enforced in the constrained formulation
}
\begin{align}
\min_{\Gamma, \theta} \ & \text{TC}(\Gamma)  \label{general_machine_teaching_con1:obj}\\
\mbox{s.t. } & \text{TR}(\theta) \leq \epsilon \revtwoCG{,} \label{general_machine_teaching_con1:TR}\\
& \theta = \text{L} (\Gamma) \revtwoCG{,} \label{general_machine_teaching_con1:ML}
\end{align}
\CGrev{which allows for either approximate or exact teaching. Alternatively, 
given a teaching budget $B$, the learning is performed via the constrained formulation}
\begin{align}
\min_{\Gamma, \theta} \ & \text{TR}(\theta) \label{general_machine_teaching_con2:obj}\\
\mbox{s.t. } & \text{TC}(\Gamma) \leq B \revtwoCG{,} \label{general_machine_teaching_con2:TC}\\
& \theta = \text{L} (\Gamma). \label{general_machine_teaching_con2:ML}
\end{align}

\CGrev{
Other variants consider multiple learners to be taught by the same teacher (i.e., common teaching set). The teacher can aim to optimize for the worst learner (minimax risk), or the average learner (Bayes risk).}
For the teaching dimension problem, the teaching cost is the cardinality of the teaching dataset, namely its $0$-norm. 
If the empirical minimization loss $\mathcal{L}$ is guiding the learning process, and $\lambda$ is the regularization weight, then teaching dimension problem can be formulated as
\begin{align}
\min_{\Gamma, \hat{\theta}} \ & \lambda \lVert \Gamma \rVert_0  \label{machine_teaching:upper}\\
\mbox{s.t. } & \lVert \hat{\theta} - \theta^* \rVert_2^2 \leq \epsilon \revtwoCG{,}\\
& \hat{\theta} \ \CGrev{\in} \ \text{argmin}_{\theta \in \Theta} \sum_{x \in X} \mathcal{L}(\theta;x) + \lambda \lVert \theta \rVert_2^2. \label{machine_teaching:lower}
\end{align}

Machine teaching approaches tailored to specific learners have also been explored in the literature. In \cite{zhu2013machine}, a method is proposed for the Bayesian learners, while \cite{patil2014optimal} focuses on Generalized Context Model learners. In \cite{mei2015using}, the bilevel optimization of machine teaching is explored to devise optimal data poisoning attacks for a broad family of learners (i.e., SVM, logistic regression, linear regression). The attacker seeks the minimum training set poisoning to attack the learned model. By using the KKT conditions of the learner's problem, the bilevel formulation is turned into a single level optimization problem.

\subsection{Empirical Model Learning}\label{sec:emp-learning}

Empirical model learning (EML) aims to integrate machine learning models in combinatorial optimization in order to support decision-making in high-complexity systems through prescriptive analytics. This goes beyond the traditional what-if approaches where a predictive model (e.g., a simulation model) is used to estimate the parameters of an optimization model. A general framework for an EML approach is provided in \cite{lombardi2017empirical} and requires the following:
\begin{itemize}
 \item \CGrev{A vector $\eta$ of $n$ decision variables $\eta_i$, with $\eta_i$ feasible over the domain $D_i$.}
 \item A mathematical encoding $h$ of the Machine Learning model.
 \item A vector $z$ of observables obtained from $h$.
 \item Logical predicates $g_j(\eta,z)$ such as mathematical programming inequalities or combinatorial restrictions in constraint programming.
 \item A cost function $f(\eta,z)$. 
\end{itemize}
EML then solves the following optimization problem
\begin{align}
\min \ & f(\eta, z)  \label{empirical_learning:obj}\\
\mbox{s.t. } & g_j(\eta,z) \quad \forall j \in J\revtwoCG{,} \label{empirical_learning:predicates}\\
& z = h(\eta) \revtwoCG{,} \label{empirical_learning:observables}\\
& \CGrev{\eta_i \in D_i \quad \forall i = 1, \dots, n.} \label{empirical_learning:decision_variables}
\end{align}
 
 The combinatorial structure of the problem is defined by \eqref{empirical_learning:obj}, \eqref{empirical_learning:predicates}, and \eqref{empirical_learning:decision_variables} while \eqref{empirical_learning:observables} embeds the empirical machine learning model in the combinatorial problem. Embedding techniques for neural networks and decision trees are presented in \cite{lombardi2017empirical} using  optimization approaches that include \CGrev{mixed integer} \CGrev{nonlinear} programming, constraint programming, and SAT Modulo Theories, and local search. 


\JSrev{\subsection{Bayesian Network Structure Learning}\label{sec:bayesian}
	Bayesian networks are a class of models that represent cause-effect relationships. These networks are learned by deriving the causal relationships from data. A Bayesian network is visually represented as a direct acyclic graph $G(N,E)$ where each of the nodes in $N$ corresponds to one variable and the edges $E$ are directional relations that indicate the cause and effect relationships among the variables. A conditional probability distribution is associated with every node/variables and along with the network structure expresses the conditional dependencies among all the variables. A main challenge in learning Bayesian networks is learning the network structure from the data. \revtwoCGblue{This is} known as the Bayesian network structure learning problem. Finding the optimal Bayesian network structure is $\mathcal{NP}$-hard \cite{chickering1996learning}. Mixed integer programming formulations of the Bayesian network structure learning have been proposed \cite{Barlett:2013} and solved by using relaxations \cite{jaakkola2010learning}, cutting planes \cite{bartlett2017integer, campos2011efficient, Cussens:2011}, and heuristics \cite{gasse2014hybrid, yuan2013learning}.}

\JSrev{The case of learning Bayesian network structures when the width of the tree is bounded by a small constant is computationally tractable \cite{nie2014advances, parviainen2014learning}. The bounded tree-width case is thus a restriction on the Bayesian network structure that limits the ability to represent exactly the underlying distribution of the data with the aim to achieve reasonable computational performance when computing the network structure. Following \cite{nie2014advances}, to formulate the Bayesian network structure learning problem with a maximum tree-width $w$, the following binary variables are defined
	$$p_{it}=\begin{cases}
	1 &\mbox{ if $P_{it}$ is the parent set of node $i$}\\
	0 &\mbox{ otherwise}
	\end{cases}
	$$
	where $i\in N$ and $P_{it}$ is a parent set for node $i$. For each node $i$, the collection of parent sets is denoted as $P_i$ and is assumed to be available (i.e.\revtwoCGblue{,} enumerated beforehand). Thus $P_{it} \in P_i$ with $t= 1,\dots,r_i$, and $r_i=|P_i|$ where $P_i\subset N$. Additional auxiliary variables $z_i\in [0,|N|]$, $v_i\in [0,|N|]$ where $|N|$ denotes the number of nodes in $N$, and $y_{ij}\in \{0,1\}$ are introduced to enforce the tree-width and directed acyclic graph conditions. The problem is formulated as
	\begin{align}
	\max & \sum_{i\in N}\sum_{t=1}^{r_i}p_{it}s_i(P_{it}) \label{bayesian_eq1}\\
	\mbox{s.t.} & \sum_{j\in N}y_{ij}\leq w, \quad \forall i \in N,\label{bayesian_eq2}\\
	& (|N|+1)y_{ij} \leq |N| + z_j - z_i \quad \forall i,j\in N,\label{bayesian_eq3}\\
	&y_{ij}+y_{ik} - y_{jk}-y_{kj} \leq 1 \quad \forall i,j,k \in N,\label{bayesian_eq4}\\
	& \sum_{t=1}^{r_i}p_{it}=1 \quad \forall i \in N,\label{bayesian_eq5}\\
	& (|N|+1)p_{it}\leq |N| + v_j - v_i \quad \forall i \in N, \ \forall t = 1,\dots, r_i, \ \forall j \in P_{it},\label{bayesian_eq6}\\
	& p_{it}\leq y_{ij} + y_{ji} \quad \forall i \in N, \ \forall t = 1,\dots, r_i, \ \forall j \in P_{it},\label{bayesian_eq7}\\
	& p_{it}\leq y_{jk} + y_{kj} \quad \forall i \in N, \ \forall t =1,\dots, r_i, \ \forall j,k \in P_{it},\label{bayesian_eq8}\\
	&z_i \in [0,|N|], \ v_i \in[0,|N|], \ y_{ij}\in\{0,1\}, \ p_{it}\in\{0,1\} \quad \forall i,j\in N, \ \forall t= 1,\dots,r_i.\label{bayesian_eq9}
	\end{align}
	The objective function \eqref{bayesian_eq1} maximizes the score of the \revtwoCG{acyclic} graph where $s_i()$ is a score function that can be efficiently computed for every node $i\in N$ \cite{campos2011efficient}. Constraints \eqref{bayesian_eq2}--\eqref{bayesian_eq4} enforce a maximum tree-width $w$ while constraints \eqref{bayesian_eq5}--\eqref{bayesian_eq6} enforce the directed acyclic graph condition. Constraints \eqref{bayesian_eq7}--\eqref{bayesian_eq8} enforce the relationship between the $p$ and $y$ variables and finally constraints \eqref{bayesian_eq9} set the variable bounds and binary conditions. Another formulation for the bounded tree-width problem has been proposed in \cite{parviainen2014learning} and includes an exponential number of constraints which are separated in a branch-and-cut framework. Both formulations however become computationally demanding as the number of features in the data set grows and with an increase in the tree-width limit. Several search heuristics have also been proposed as  solution approaches \cite{nie2014advances, nie2015learning, scanagatta2016learning}.
}

\section{Conclusions}\label{sec:conclusion}

Mathematical programming constitutes a fundamental aspect of many machine learning models where the training of these models is a large scale optimization problem. This paper surveyed a wide range of machine learning models namely regression, classification, clustering, and deep learning as well as the new emerging paradigms of machine teaching and empirical model learning. The important mathematical optimization models for expressing these machine learning models are presented and discussed. Exploiting the large scale optimization formulations and devising model specific solution approaches is an important line of research particularly benefiting from the maturity of commercial optimization software \JSrev{to solve the problems to optimality or to devise effective heuristics. \revtwoCG{However, as highlighted in \cite{liang2017fisher, poggio2017theory},} providing quantitative performance bounds remains an open problem.} The nonlinearity of the models, the associated uncertainty of the data, as well as the scale of the problems represent some of the very important and compelling challenges to the mathematical optimization community. Furthermore, bilevel formulations play a big role in adversarial learning  \cite{HamN18}, including adversarial training, data poisoning and neural network robustness.

\js{Based on this survey, we summarize the distinctive features and the potential open machine learning problems that may benefit from the advances in computational optimization.
\begin{itemize}
	\item \textbf{Regression.} The typical approaches to avoid overfitting and to handle uncertainty in the data include shrinkage methods and dimension reduction. These approaches can all be posed as mathematical programming models. General non-convex regularization to enforce sparsity without incurring shrinkage and bias (such as in lasso and ridge regularization) remain computationally challenging to solve to optimality. Investigating tighter relaxations and exact solution approaches continue to be an active line of research \cite{atamturk2019rank}.
	\item \textbf{Classification.} Classification problems can also be naturally formulated as optimization problems. Support vector machines in particular have been well studied in the optimization literature. Similar to regression, classifier sparsity is one important approach to avoid overfitting. Additionally, exploiting the kernel tricks is key as nonlinear separators are obtained without additional complexity. However, when posed as an optimization problem, it is still unclear how to exploit kernel tricks in sparse SVM optimization models. Another advantage to express machine learning problems as optimization problems and in particular classification problems is to account for inaccuracies in the data. Handling data uncertainty is a \revtwoCGblue{deeply explored} field in the optimization literature and several practical approaches have been presented to handle uncertainty through robust and stochastic optimization. Such advances in the optimization literature are currently being investigated to improve over the standard approaches \cite{bertismas2019robustclassification}.
	\item \textbf{Clustering.} Clustering problems are in general formulated as MINLPs that are hard to solve to optimality. The challenges include handling the non-convexity as well as the large scale instances which is a challenge even for linear variants such as the capacitated centred clustering (formulated as a binary linear model). Especially for large-scale instances, heuristics are typically devised. Exact approaches for clustering received less attention in the literature.	
	\item \textbf{DNNs architectures as MIPs}. The advantage of mathematical programming approaches to model DNNs has only been showcased for relatively small size data sets due to the scale of the underlying optimization model. Furthermore, expressing misclassification conditions for adversarial examples in a non-restrictive manner, and handling the uncertainty in the training data are open problems in this context. 
	\item \textbf{Adversarial learning and adversarial robustness.} 
	Optimization models for the search for adversarial examples are important to identify and subsequently protect against novel sets of attacks. The complexity of the mathematical models in this context is highly dependent on the the classifier function. Untargeted attacks received less attention in the literature, and the mathematical programming formulation \eqref{adv_dist_bi:obj}--\eqref{adv_dist_bi:fin} has been introduced in \cref{sec:untargeted-attacks}. Furthermore, designing models robust to adversarial attacks is a two-player game, which can be cast as a bilevel optimization problem. The loss function adopted by the learner is one main   complexity for the resulting mathematical model and solution approaches remain to be investigated.
	\item \textbf{Data poisoning}: Similar to adversarial robustness, defending against the poisoning of the training data is a two-player game.  
	The case of online data retrieval is especially challenging for gradient-based algorithms as the KKT conditions do not hold.
	\item \textbf{Activation ensembles.} 
	Activation ensembles seek a trade-off between the classifier accuracy and computational feasibility of training with a mathematical programming approach. Adopting activation ensembles to train large DNNs have not been investigated yet.
	\item \textbf{Machine teaching.} Posed as a bilevel optimization problem, one of the challenges in machine teaching is to devise computationally tractable single-level formulations that model the learner, the teaching risk, and the teaching cost. Machine teaching also generalizes a number of two-player games that are important in practice including data poisoning and adversarial training.
	\item \textbf{Empirical model learning.} This emerging paradigm can be seen as the bridge combining machine learning for parameter estimation and operations research for optimization. As such, theoretical and practical challenges remain to be investigated to propose prescriptive analytics models jointly combining learning and optimization in practical applications.
\end{itemize}
}

While this survey does not discuss numerical optimization techniques since they were recently reviewed in \cite{bottou2018optimization, Cur17, wright2017optimization}, we note the fundamental role of the stochastic gradient algorithm \cite{robbins1985stochastic} \cg{and its variants} on large scale machine learning. We also highlight the potential impact of machine learning on advancing the solution approaches of mathematical programming \cite{FISCHETTI2018, fischetti2018learning}. 

\JSrev{
This survey has also focused on the learning process (loss minimization), however we note that challenging optimization problems also appear in the inference process, i.e., energy minimization (see \cite{lecun_energy} for a comprehensive survey). In the inference step, the best output $y^*$ is chosen from among all possible outputs given a certain input $x$ such that an ``energy function'' is minimized. The energy function provides a measure of the goodness of a particular configuration of the input and output variables. Energy optimization constitute a common framework for machine learning where the training of a model aims at finding the optimal energy function.} 

\CGrev{A key part of most machine learning approaches is the choice of the {hyperparameters} of the learning model.
	The Hyperparameter Optimization (HPO) is usually driven by the data scientist's experience and the characteristics of the dataset and typically follows heuristic rules or cross-validation approaches. Alternatively, \revtwoCG{the HPO problem} can be modeled as a box-constrained mathematical optimization problem \cite{diaz2017effective}, or as a bilevel optimization problem as discussed in \cite{franceschi2018bilevel, klatzer2015continuous, moore2011model}, which provides theoretical convergence guarantees in addition to computational advantage. Automated approaches for HPO are also an active area of research in Machine Learning \cite{bergstra2012random, elsken2019neural, wistuba2019survey}.}
	
 Finally, since the recent widespread of machine learning to several research disciplines and in the mainstream industry can be largely attributed to the availability of data and the relatively easy to use libraries, we summarize in the online supplement the resources that may be of value for research.

%
%

\section*{Acknowledgement} We are very grateful to four anonymous referees for their valuable feedback and comments that helped improve the content and presentation of the paper. Joe Naoum-Sawaya was supported by NSERC Discovery Grant RGPIN-2017-03962 and Bissan Ghaddar was supported by NSERC Discovery Grant RGPIN-2017-04185.

\bibliographystyle{plain}

\bibliography{opt_ML,datasets,learning_frameworks}

\begin{thebibliography}{100}

\bibitem{agos2014}
Forest Agostinelli, Matthew Hoffman, Peter Sadowski, and Pierre Baldi.
\newblock Learning activation functions to improve deep neural networks.
\newblock Technical report, arXiv preprint 1412.6830, 2014.

\bibitem{alam2018kernel}
Md~Ashad Alam, Hui-Yi Lin, Hong-Wen Deng, Vince~D Calhoun, and Yu-Ping Wang.
\newblock A kernel machine method for detecting higher order interactions in
  multimodal datasets: Application to schizophrenia.
\newblock {\em Journal of Neuroscience Methods}, 309:161--174, 2018.

\bibitem{Aloise2012}
Daniel Aloise, Pierre Hansen, and Leo Liberti.
\newblock An improved column generation algorithm for minimum sum-of-squares
  clustering.
\newblock {\em Mathematical Programming}, 131(1):195--220, 2012.

\bibitem{amaldi2013distance}
Edoardo Amaldi and Stefano Coniglio.
\newblock A distance-based point-reassignment heuristic for the k-hyperplane
  clustering problem.
\newblock {\em European Journal of Operational Research}, 227(1):22--29, 2013.

\bibitem{amaldi2016discrete}
Edoardo Amaldi, Stefano Coniglio, and Leonardo Taccari.
\newblock Discrete optimization methods to fit piecewise affine models to data
  points.
\newblock {\em Computers \& Operations Research}, 75:214--230, 2016.

\bibitem{aoki2014cluster}
Yasunori Aoki, Ken Hayami, Hans~De Sterck, and Akihiko Konagaya.
\newblock Cluster {N}ewton method for sampling multiple solutions of
  underdetermined inverse problems: application to a parameter identification
  problem in pharmacokinetics.
\newblock {\em SIAM Journal on Scientific Computing}, 36(1):14--44, 2014.

\bibitem{atamturk2019rank}
Alper Atamturk and Andres Gomez.
\newblock Rank-one convexification for sparse regression.
\newblock Technical report, arXiv preprint 1901.10334, 2019.

\bibitem{aytug2015feature}
Haldun Aytug.
\newblock Feature selection for support vector machines using generalized
  {B}enders decomposition.
\newblock {\em European Journal of Operational Research}, 244(1):210--218,
  2015.

\bibitem{6982798}
M.~Azad and M.~Moshkov.
\newblock Minimization of decision tree depth for multi-label decision tables.
\newblock In {\em Proceedings of the IEEE International Conference on Granular
  Computing}, pages 7--12, 2014.

\bibitem{7321423}
M.~Azad and M.~Moshkov.
\newblock Classification and optimization of decision trees for inconsistent
  decision tables represented as {MVD} tables.
\newblock In {\em Proceedings of the Federated Conference on Computer Science
  and Information Systems}, pages 31--38, 2015.

\bibitem{AZAD2014368}
Mohammad Azad and Mikhail Moshkov.
\newblock Minimization of decision tree average depth for decision tables with
  many-valued decisions.
\newblock {\em Procedia Computer Science}, 35:368--377, 2014.

\bibitem{AZAD2017910}
Mohammad Azad and Mikhail Moshkov.
\newblock Multi-stage optimization of decision and inhibitory trees for
  decision tables with many-valued decisions.
\newblock {\em European Journal of Operational Research}, 263(3):910--921,
  2017.

\bibitem{bagirov2006new}
Adil~M Bagirov and John Yearwood.
\newblock A new nonsmooth optimization algorithm for minimum sum-of-squares
  clustering problems.
\newblock {\em European Journal of Operational Research}, 170(2):578--596,
  2006.

\bibitem{Barlett:2013}
Mark Barlett and James Cussens.
\newblock Advances in {B}ayesian network learning using integer programming.
\newblock In {\em Proceedings of the Conference on Uncertainty in Artificial
  Intelligence}, pages 182--191, 2013.

\bibitem{barreno2010security}
Marco Barreno, Blaine Nelson, Anthony~D. Joseph, and J.~Doug Tygar.
\newblock The security of machine learning.
\newblock {\em Machine Learning}, 81(2):121--148, 2010.

\bibitem{bartlett2017integer}
Mark Bartlett and James Cussens.
\newblock Integer linear programming for the {B}ayesian network structure
  learning problem.
\newblock {\em Artificial Intelligence}, 244:258--271, 2017.

\bibitem{bastani2016measuring}
Osbert Bastani, Yani Ioannou, Leonidas Lampropoulos, Dimitrios Vytiniotis,
  Aditya Nori, and Antonio Criminisi.
\newblock Measuring neural net robustness with constraints.
\newblock In D.~D. Lee, M.~Sugiyama, U.~V. Luxburg, I.~Guyon, and R.~Garnett,
  editors, {\em Advances in Neural Information Processing Systems}, pages
  2613--2621. Curran Associates, Inc., 2016.

\bibitem{BAUMANN20191041}
Philipp Baumann, D.~S. Hochbaum, and Y.~T. Yang.
\newblock A comparative study of the leading machine learning techniques and
  two new optimization algorithms.
\newblock {\em European Journal of Operational Research}, 272(3):1041--1057,
  2019.

\bibitem{belhumeur1996eigenfaces}
Peter~N Belhumeur, Jo{\~a}o~P Hespanha, and David~J Kriegman.
\newblock Eigenfaces vs. fisherfaces: Recognition using class specific linear
  projection.
\newblock {\em IEEE Transactions on Pattern Analysis \& Machine Intelligence},
  (7):711--720, 1997.

\bibitem{benati2014mixed}
Stefano Benati and Sergio Garc{\'\i}a.
\newblock A mixed integer linear model for clustering with variable selection.
\newblock {\em Computers \& Operations Research}, 43:280--285, 2014.

\bibitem{bengio2018machine}
Yoshua Bengio, Andrea Lodi, and Antoine Prouvost.
\newblock Machine learning for combinatorial optimization: a methodological
  tour d'horizon.
\newblock Technical report, arXiv preprint 1811.06128, 2018.

\bibitem{bennett1992decision}
Kristin~P Bennett.
\newblock Decision tree construction via linear programming.
\newblock Technical report, Center for Parallel Optimization, Computer Sciences
  Department, University of Wisconsin, 1992.

\bibitem{bennett1996optimal}
Kristin~P. Bennett and J.~Blue.
\newblock Optimal decision trees.
\newblock Technical report, Rensselaer Polytechnic Institute, 1996.

\bibitem{bennett1992robust}
Kristin~P Bennett and Olvi~L Mangasarian.
\newblock Robust linear programming discrimination of two linearly inseparable
  sets.
\newblock {\em Optimization Methods and Software}, 1(1):23--34, 1992.

\bibitem{bennett2006interplay}
Kristin~P Bennett and Emilio Parrado-Hern{\'a}ndez.
\newblock The interplay of optimization and machine learning research.
\newblock {\em Journal of Machine Learning Research}, 7:1265--1281, 2006.

\bibitem{bergstra2012random}
James Bergstra and Yoshua Bengio.
\newblock Random search for hyper-parameter optimization.
\newblock {\em Journal of Machine Learning Research}, 13(Feb):281--305, 2012.

\bibitem{BERTSIMAS2018931}
Dimitris Bertsimas and Martin~S. Copenhaver.
\newblock Characterization of the equivalence of robustification and
  regularization in linear and matrix regression.
\newblock {\em European Journal of Operational Research}, 270(3):931--942,
  2018.

\bibitem{bertsimas2017optimal}
Dimitris Bertsimas and Jack Dunn.
\newblock Optimal classification trees.
\newblock {\em Machine Learning}, 106(7):1039--1082, 2017.

\bibitem{bertismas2019robustclassification}
Dimitris Bertsimas, Jack Dunn, Colin Pawlowski, and Ying~Daisy Zhuo.
\newblock Robust classification.
\newblock {\em {INFORMS} Journal on Optimization}, 1(1):2--34, 2019.

\bibitem{bertsimas2014predictive}
Dimitris Bertsimas and Nathan Kallus.
\newblock From predictive to prescriptive analytics.
\newblock {\em Management Science}, 66(3):1025--1044, 2020.

\bibitem{bertsimasORforum}
Dimitris Bertsimas and Angela King.
\newblock {OR} forum--{A}n algorithmic approach to linear regression.
\newblock {\em Operations Research}, 64(1):2--16, 2016.

\bibitem{bertsimas2016best}
Dimitris Bertsimas, Angela King, and Rahul Mazumder.
\newblock Best subset selection via a modern optimization lens.
\newblock {\em The Annals of Statistics}, 44(2):813--852, 2016.

\bibitem{bertsimas2007CR}
Dimitris Bertsimas and Romy Shioda.
\newblock Classification and regression via integer optimization.
\newblock {\em Operations Research}, 55(2):252--271, 2007.

\bibitem{bertsimas2017sparseregr}
Dimitris Bertsimas, Bart Van~Parys, et~al.
\newblock Sparse high-dimensional regression: Exact scalable algorithms and
  phase transitions.
\newblock {\em The Annals of Statistics}, 48(1):300--323, 2020.

\bibitem{biggio2010multiple}
Battista Biggio, Giorgio Fumera, and Fabio Roli.
\newblock Multiple classifier systems under attack.
\newblock In {\em Proceedings of the International Workshop on Multiple
  Classifier Systems}, pages 74--83, 2010.

\bibitem{biggio2012poisoning}
Battista Biggio, Blaine Nelson, and Pavel Laskov.
\newblock Poisoning attacks against support vector machines.
\newblock In {\em Proceedings of the International Conference on Machine
  Learning}, pages 1467--1474, 2012.

\bibitem{blanco2018locating}
V{\'\i}ctor Blanco, Justo Puerto, and Rom{\'a}n Salmer{\'o}n.
\newblock Locating hyperplanes to fitting set of points: A general framework.
\newblock {\em Computers \& Operations Research}, 95:172--193, 2018.

\bibitem{blanquero2018optimal}
Rafael Blanquero, Emilio Carrizosa, Cristina Molero-R{\i}o, and Dolores~Romero
  Morales.
\newblock Optimal randomized classification trees.
\newblock Technical report, 2018.

\bibitem{blanquero2018sparsity}
Rafael Blanquero, Emilio Carrizosa, Cristina Molero-R{\'\i}o, and
  Dolores~Romero Morales.
\newblock Sparsity in optimal randomized classification trees.
\newblock {\em European Journal of Operational Research}, 284(1):255--272,
  2020.

\bibitem{Bonami2015}
Pierre Bonami, Andrea Lodi, Andrea Tramontani, and Sven Wiese.
\newblock On mathematical programming with indicator constraints.
\newblock {\em Mathematical Programming}, 151(1):191--223, 2015.

\bibitem{BonLZ18}
Pierre Bonami, Andrea Lodi, and Giulia Zarpellon.
\newblock Learning a classification of mixed-integer quadratic programming
  problems.
\newblock In {\em Proceedings of the International Conference on the
  Integration of Constraint Programming, Artificial Intelligence, and
  Operations Research}, pages 595--604, 2018.

\bibitem{boct2011optimization}
Radu~Ioan Bo{\c{t}} and Nicole Lorenz.
\newblock Optimization problems in statistical learning: Duality and optimality
  conditions.
\newblock {\em European Journal of Operational Research}, 213(2):395--404,
  2011.

\bibitem{bottou2018optimization}
L{\'e}on Bottou, Frank~E. Curtis, and Jorge Nocedal.
\newblock Optimization methods for large-scale machine learning.
\newblock {\em SIAM Review}, 60(2):223--311, 2018.

\bibitem{bradley2000massive}
Paul Bradley and Olvi Mangasarian.
\newblock Massive data discrimination via linear support vector machines.
\newblock {\em Optimization Methods and Software}, 13(1):1--10, 2000.

\bibitem{breiman1984classification}
L~Breiman, J~Friedman, R~Olshen, and C~Stone.
\newblock {\em Classification and Regression Trees}.
\newblock {Chapman and Hall/CRC, London}, 1984.

\bibitem{bruckner2012static}
Michael Br{\"u}ckner, Christian Kanzow, and Tobias Scheffer.
\newblock Static prediction games for adversarial learning problems.
\newblock {\em Journal of Machine Learning Research}, 13:2617--2654, 2012.

\bibitem{bruckner2011stackelberg}
Michael Br{\"u}ckner and Tobias Scheffer.
\newblock Stackelberg games for adversarial prediction problems.
\newblock In {\em Proceedings of the International Conference on Knowledge
  Discovery and Data Mining}, pages 547--555, 2011.

\bibitem{BunTTKK18}
Rudy~R Bunel, Ilker Turkaslan, Philip Torr, Pushmeet Kohli, and Pawan~K
  Mudigonda.
\newblock A unified view of piecewise linear neural network verification.
\newblock In S.~Bengio, H.~Wallach, H.~Larochelle, K.~Grauman, N.~Cesa-Bianchi,
  and R.~Garnett, editors, {\em Advances in Neural Information Processing
  Systems}, pages 4790--4799. Curran Associates, Inc., 2018.

\bibitem{byrd1995limited}
Richard~H Byrd, Peihuang Lu, Jorge Nocedal, and Ciyou Zhu.
\newblock A limited memory algorithm for bound constrained optimization.
\newblock {\em SIAM Journal on Scientific Computing}, 16(5):1190--1208, 1995.

\bibitem{cafieri2014reformulation}
Sonia Cafieri, Alberto Costa, and Pierre Hansen.
\newblock Reformulation of a model for hierarchical divisive graph modularity
  maximization.
\newblock {\em Annals of Operations Research}, 222(1):213--226, 2014.

\bibitem{CAFIERI201465}
Sonia Cafieri, Pierre Hansen, and Leo Liberti.
\newblock Improving heuristics for network modularity maximization using an
  exact algorithm.
\newblock {\em Discrete Applied Mathematics}, 163:65--72, 2014.

\bibitem{campos2011efficient}
Cassio P~de Campos and Qiang Ji.
\newblock Efficient structure learning of {B}ayesian networks using
  constraints.
\newblock {\em Journal of Machine Learning Research}, 12:663--689, 2011.

\bibitem{carlini2017towards}
Nicholas Carlini and David Wagner.
\newblock Towards evaluating the robustness of neural networks.
\newblock In {\em Proceedings of the IEEE Symposium on Security and Privacy},
  pages 39--57, 2017.

\bibitem{CARRIZOSA2014151}
Emilio Carrizosa and Vanesa Guerrero.
\newblock Biobjective sparse principal component analysis.
\newblock {\em Journal of Multivariate Analysis}, 132:151--159, 2014.

\bibitem{CARRIZOSA2014349}
Emilio Carrizosa and Vanesa Guerrero.
\newblock rs-{S}parse principal component analysis: A mixed integer nonlinear
  programming approach with {VNS}.
\newblock {\em Computers \& Operations Research}, 52:349--354, 2014.

\bibitem{CarrIJOC}
Emilio Carrizosa, Belén Martín-Barragán, and Dolores~Romero Morales.
\newblock Binarized support vector machines.
\newblock {\em {INFORMS} Journal on Computing}, 22(1):154--167, 2010.

\bibitem{CARRIZOSA2011260}
Emilio Carrizosa, Belén Martín-Barragán, and Dolores~Romero Morales.
\newblock Detecting relevant variables and interactions in supervised
  classification.
\newblock {\em European Journal of Operational Research}, 213(1):260--269,
  2011.

\bibitem{CARRIZOSA2013356}
Emilio Carrizosa, Nenad Mladenović, and Raca Todosijević.
\newblock Variable neighborhood search for minimum sum-of-squares clustering on
  networks.
\newblock {\em European Journal of Operational Research}, 230(2):356--363,
  2013.

\bibitem{carrizosa2013supervised}
Emilio Carrizosa and Dolores~Romero Morales.
\newblock Supervised classification and mathematical optimization.
\newblock {\em Computers \& Operations Research}, 40(1):150--165, 2013.

\bibitem{chan2007direct}
Antoni~B. Chan, Nuno Vasconcelos, and Gert R.~G. Lanckriet.
\newblock Direct convex relaxations of sparse {SVM}.
\newblock In {\em Proceedings of the International Conference on Machine
  Learning}, pages 145--153, 2007.

\bibitem{chatterjee2015regression}
Samprit Chatterjee and Ali~S Hadi.
\newblock {\em Regression analysis by example}.
\newblock John Wiley \& Sons, New York, 2015.

\bibitem{chaves2010clustering}
Antonio~Augusto Chaves and Luiz Antonio~Nogueira Lorena.
\newblock Clustering search algorithm for the capacitated centered clustering
  problem.
\newblock {\em Computers \& Operations Research}, 37(3):552--558, 2010.

\bibitem{chen2013complete}
Xiaobo Chen, Jian Yang, David Zhang, and Jun Liang.
\newblock Complete large margin linear discriminant analysis using mathematical
  programming approach.
\newblock {\em Pattern Recognition}, 46(6):1579--1594, 2013.

\bibitem{chen1995nonlinear}
Yang Chen and Michael Florian.
\newblock The nonlinear bilevel programming problem: Formulations, regularity
  and optimality conditions.
\newblock {\em Optimization}, 32(3):193--209, 1995.

\bibitem{cheng2017maximum}
Chih-Hong Cheng, Georg N{\"u}hrenberg, and Harald Ruess.
\newblock Maximum resilience of artificial neural networks.
\newblock In Deepak D'Souza and K.~Narayan~Kumar, editors, {\em Automated
  Technology for Verification and Analysis}, pages 251--268, Cham, 2017.
  Springer International Publishing.

\bibitem{chickering1996learning}
David~Maxwell Chickering.
\newblock Learning {B}ayesian networks is {NP}-complete.
\newblock In {\em Learning from Data}, pages 121--130. Springer, 1996.

\bibitem{CHIKALOV2018689}
Igor Chikalov, Shahid Hussain, and Mikhail Moshkov.
\newblock Bi-criteria optimization of decision trees with applications to data
  analysis.
\newblock {\em European Journal of Operational Research}, 266(2):689--701,
  2018.

\bibitem{chouldechova2015generalized}
Alexandra Chouldechova and Trevor Hastie.
\newblock Generalized additive model selection.
\newblock Technical report, arXiv preprint 1506.03850, 2015.

\bibitem{chu2007support}
Wei Chu and S~Sathiya Keerthi.
\newblock Support vector ordinal regression.
\newblock {\em Neural Computation}, 19(3):792--815, 2007.

\bibitem{claassen2007application}
GDH Claassen and Th~HB Hendriks.
\newblock An application of special ordered sets to a periodic milk collection
  problem.
\newblock {\em European Journal of Operational Research}, 180(2):754--769,
  2007.

\bibitem{corne2012synergies}
David Corne, Clarisse Dhaenens, and Laetitia Jourdan.
\newblock Synergies between operations research and data mining: The emerging
  use of multi-objective approaches.
\newblock {\em European Journal of Operational Research}, 221(3):469--479,
  2012.

\bibitem{corrente2013robust}
Salvatore Corrente, Salvatore Greco, Mi{\l}osz Kadzi{\'n}ski, and Roman
  S{\l}owi{\'n}ski.
\newblock Robust ordinal regression in preference learning and ranking.
\newblock {\em Machine Learning}, 93(2-3):381--422, 2013.

\bibitem{Cortes1995}
Corinna Cortes and Vladimir Vapnik.
\newblock Support-vector networks.
\newblock {\em Machine Learning}, 20(3):273--297, Sep 1995.

\bibitem{courbariaux2015binaryconnect}
Matthieu Courbariaux, Yoshua Bengio, and Jean-Pierre David.
\newblock Binaryconnect: Training deep neural networks with binary weights
  during propagations.
\newblock In C.~Cortes, N.~D. Lawrence, D.~D. Lee, M.~Sugiyama, and R.~Garnett,
  editors, {\em Advances in Neural Information Processing Systems}, pages
  3123--3131. Curran Associates, Inc., 2015.

\bibitem{cox1989heuristic}
Louis~Anthony Cox, Yuping Qiu, and Warren Kuehner.
\newblock Heuristic least-cost computation of discrete classification functions
  with uncertain argument values.
\newblock {\em Annals of Operations Research}, 21(1):1--29, 1989.

\bibitem{cunningham2015linear}
John~P Cunningham and Zoubin Ghahramani.
\newblock Linear dimensionality reduction: survey, insights, and
  generalizations.
\newblock {\em The Journal of Machine Learning Research}, 16(1):2859--2900,
  2015.

\bibitem{Cur17}
Frank~E. Curtis and Katya Scheinberg.
\newblock Optimization methods for supervised machine learning: From linear
  models to deep learning.
\newblock In {\em Leading Developments from {INFORMS} Communities}, pages
  89--114. INFORMS, 2017.

\bibitem{Cussens:2011}
James Cussens.
\newblock {B}ayesian network learning with cutting planes.
\newblock In {\em Proceedings of the Conference on Uncertainty in Artificial
  Intelligence}, pages 153--160, 2011.

\bibitem{cybenko1989approximation}
George Cybenko.
\newblock Approximation by superpositions of a sigmoidal function.
\newblock {\em Mathematics of Control, Signals and Systems}, 2(4):303--314,
  1989.

\bibitem{d2015mathematical}
Claudia D'Ambrosio, Andrea Lodi, Sven Wiese, and Cristiana Bragalli.
\newblock Mathematical programming techniques in water network optimization.
\newblock {\em European Journal of Operational Research}, 243(3):774--788,
  2015.

\bibitem{DEFARIAS2008234}
I.R. de~Farias, M.~Zhao, and H.~Zhao.
\newblock A special ordered set approach for optimizing a discontinuous
  separable piecewise linear function.
\newblock {\em Operations Research Letters}, 36(2):234--238, 2008.

\bibitem{dekel2010learning}
Ofer Dekel, Ohad Shamir, and Lin Xiao.
\newblock Learning to classify with missing and corrupted features.
\newblock {\em Machine Learning}, 81(2):149--178, 2010.

\bibitem{diaz2017effective}
Gonzalo~I. Diaz, Achille Fokoue-Nkoutche, Giacomo Nannicini, and Horst
  Samulowitz.
\newblock An effective algorithm for hyperparameter optimization of neural
  networks.
\newblock {\em IBM Journal of Research and Development}, 61(4):9--1, 2017.

\bibitem{diaz2004continuous}
JM~D{\i}az-B{\'a}nez, Juan~A Mesa, and Anita Sch{\"o}bel.
\newblock Continuous location of dimensional structures.
\newblock {\em European Journal of Operational Research}, 152(1):22--44, 2004.

\bibitem{ding2004k}
Chris Ding and Xiaofeng He.
\newblock K-means clustering via principal component analysis.
\newblock In {\em Proceedings of the International Conference on Machine
  Learning}, page~29, 2004.

\bibitem{doshi2017towards}
Finale Doshi-Velez and Been Kim.
\newblock Towards a rigorous science of interpretable machine learning.
\newblock Technical report, arXiv preprint 1702.08608, 2017.

\bibitem{dreiseitl2002logistic}
Stephan Dreiseitl and Lucila Ohno-Machado.
\newblock Logistic regression and artificial neural network classification
  models: a methodology review.
\newblock {\em Journal of Biomedical Informatics}, 35(5-6):352--359, 2002.

\bibitem{DUNBAR2010470}
Michelle Dunbar, John~M. Murray, Lucette~A. Cysique, Bruce~J. Brew, and
  Vaithilingam Jeyakumar.
\newblock Simultaneous classification and feature selection via convex
  quadratic programming with application to {HIV}-associated neurocognitive
  disorder assessment.
\newblock {\em European Journal of Operational Research}, 206(2):470--478,
  2010.

\bibitem{edgeworth1887observations}
Francis~Y Edgeworth.
\newblock On observations relating to several quantities.
\newblock {\em Hermathena}, 6(13):279--285, 1887.

\bibitem{edmunds1992algorithm}
Thomas~A Edmunds and Jonathan~F Bard.
\newblock An algorithm for the mixed-integer nonlinear bilevel programming
  problem.
\newblock {\em Annals of Operations Research}, 34(1):149--162, 1992.

\bibitem{efron2004least}
Bradley Efron, Trevor Hastie, Iain Johnstone, Robert Tibshirani, et~al.
\newblock Least angle regression.
\newblock {\em The Annals of Statistics}, 32(2):407--499, 2004.

\bibitem{elsken2019neural}
Thomas Elsken, Jan~Hendrik Metzen, and Frank Hutter.
\newblock Neural architecture search: A survey.
\newblock {\em Journal of Machine Learning Research}, 20(55):1--21, 2019.

\bibitem{fanghanel2009bilevel}
Diana Fangh{\"a}nel and Stephan Dempe.
\newblock Bilevel programming with discrete lower level problems.
\newblock {\em Optimization}, 58(8):1029--1047, 2009.

\bibitem{ferrari2003clustering}
Giancarlo Ferrari-Trecate, Marco Muselli, Diego Liberati, and Manfred Morari.
\newblock A clustering technique for the identification of piecewise affine
  systems.
\newblock {\em Automatica}, 39(2):205--217, 2003.

\bibitem{FISCHETTI2018}
Martina Fischetti and Marco Fraccaro.
\newblock Machine learning meets mathematical optimization to predict the
  optimal production of offshore wind parks.
\newblock {\em Computers \& Operations Research}, 106:289--297, 2019.

\bibitem{fischetti2018learning}
Martina Fischetti, Andrea Lodi, and Giulia Zarpellon.
\newblock Learning {MILP} resolution outcomes before reaching time-limit.
\newblock In {\em Proceedings of the International Conference on Integration of
  Constraint Programming, Artificial Intelligence, and Operations Research},
  pages 275--291, 2019.

\bibitem{fischetti2018deep}
Matteo Fischetti and Jason Jo.
\newblock Deep neural networks and mixed integer linear optimization.
\newblock {\em Constraints}, 23(3):296--309, 2018.

\bibitem{franceschi2018bilevel}
Luca Franceschi, Paolo Frasconi, Saverio Salzo, Riccardo Grazzi, and
  Massimilano Pontil.
\newblock Bilevel programming for hyperparameter optimization and
  meta-learning.
\newblock Technical report, arXiv preprint 1806.04910, 2018.

\bibitem{friedman2001elements}
Jerome Friedman, Trevor Hastie, and Robert Tibshirani.
\newblock {\em The Elements of Statistical Learning}, volume~1.
\newblock Springer Series in Statistics New York, NY, USA, 2001.

\bibitem{fukunaga2013introduction}
Keinosuke Fukunaga.
\newblock {\em Introduction to Statistical Pattern Recognition}.
\newblock Elsevier, 2013.

\bibitem{ganesh2007cloves}
K~Ganesh and TT~Narendran.
\newblock Cloves: A cluster-and-search heuristic to solve the vehicle routing
  problem with delivery and pick-up.
\newblock {\em European Journal of Operational Research}, 178(3):699--717,
  2007.

\bibitem{gasse2014hybrid}
Maxime Gasse, Alex Aussem, and Haytham Elghazel.
\newblock A hybrid algorithm for {B}ayesian network structure learning with
  application to multi-label learning.
\newblock {\em Expert Systems with Applications}, 41(15):6755--6772, 2014.

\bibitem{gaudioso2017lagrangian}
Manlio Gaudioso, Enrico Gorgone, Martine Labb{\'e}, and Antonio~M
  Rodr{\'\i}guez-Ch{\'\i}a.
\newblock {L}agrangian relaxation for {SVM} feature selection.
\newblock {\em Computers \& Operations Research}, 87:137--145, 2017.

\bibitem{gaudreau2015improvements}
Philippe Gaudreau, Ken Hayami, Yasunori Aoki, Hassan Safouhi, and Akihiko
  Konagaya.
\newblock Improvements to the cluster {N}ewton method for underdetermined
  inverse problems.
\newblock {\em Journal of Computational and Applied Mathematics}, 283:122--141,
  2015.

\bibitem{ghaddar2018high}
Bissan Ghaddar and Joe Naoum-Sawaya.
\newblock High dimensional data classification and feature selection using
  support vector machines.
\newblock {\em European Journal of Operational Research}, 265(3):993--1004,
  2018.

\bibitem{globerson2006nightmare}
Amir Globerson and Sam Roweis.
\newblock Nightmare at test time: robust learning by feature deletion.
\newblock In {\em Proceedings of the International Conference on Machine
  Learning}, pages 353--360, 2006.

\bibitem{GOLDMAN199520}
Sally Goldman and Michael Kearns.
\newblock On the complexity of teaching.
\newblock {\em Journal of Computer and System Sciences}, 50(1):20--31, 1995.

\bibitem{goodfellow2016deep}
Ian Goodfellow, Yoshua Bengio, Aaron Courville, and Yoshua Bengio.
\newblock {\em Deep Learning}, volume~1.
\newblock MIT Press Cambridge, 2016.

\bibitem{goodfellow2014generative}
Ian Goodfellow, Jean Pouget-Abadie, Mehdi Mirza, Bing Xu, David Warde-Farley,
  Sherjil Ozair, Aaron Courville, and Yoshua Bengio.
\newblock Generative adversarial nets.
\newblock In Z.~Ghahramani, M.~Welling, C.~Cortes, N.~D. Lawrence, and K.~Q.
  Weinberger, editors, {\em Advances in Neural Information Processing Systems},
  pages 2672--2680. Curran Associates, Inc., 2014.

\bibitem{goodfellow2013maxout}
Ian~J. Goodfellow, David Warde-Farley, Mehdi Mirza, Aaron Courville, and Yoshua
  Bengio.
\newblock Maxout networks.
\newblock In {\em Proceedings of the International Conference on Machine
  Learning}, pages 1319--1327, 2013.

\bibitem{grossmann2002review}
Ignacio~E. Grossmann.
\newblock Review of nonlinear mixed-integer and disjunctive programming
  techniques.
\newblock {\em Optimization and Engineering}, 3(3):227--252, 2002.

\bibitem{gu2014towards}
Shixiang Gu and Luca Rigazio.
\newblock Towards deep neural network architectures robust to adversarial
  examples.
\newblock Technical report, arXiv preprint 1412.5068, 2014.

\bibitem{gumucs2001global}
Zeynep~H G{\"u}m{\"u}{\c{s}} and Christodoulos~A Floudas.
\newblock Global optimization of nonlinear bilevel programming problems.
\newblock {\em Journal of Global Optimization}, 20(1):1--31, 2001.

\bibitem{gunluk2018optimal}
Oktay G{\"u}nl{\"u}k, Jayant Kalagnanam, Matt Menickelly, and Katya Scheinberg.
\newblock Optimal decision trees for categorical data via integer programming.
\newblock Technical report, Optimization Online, 2018.

\bibitem{guyon2002gene}
Isabelle Guyon, Jason Weston, Stephen Barnhill, and Vladimir Vapnik.
\newblock Gene selection for cancer classification using support vector
  machines.
\newblock {\em Machine Learning}, 46(1-3):389--422, 2002.

\bibitem{HamN18}
Jihun Hamm and Yung{-}Kyun Noh.
\newblock K-beam subgradient descent for minimax optimization.
\newblock Technical report, arXiv preprint 1805.11640, 2018.

\bibitem{hansen1997cluster}
Pierre Hansen and Brigitte Jaumard.
\newblock Cluster analysis and mathematical programming.
\newblock {\em Mathematical Programming}, 79(1-3):191--215, 1997.

\bibitem{har2003constraint}
Sariel Har-Peled, Dan Roth, and Dav Zimak.
\newblock Constraint classification for multiclass classification and ranking.
\newblock In S.~Becker, S.~Thrun, and K.~Obermayer, editors, {\em Advances in
  Neural Information Processing Systems}, pages 809--816. MIT Press, 2003.

\bibitem{GAMs}
Trevor Hastie and Robert Tibshirani.
\newblock Generalized additive models.
\newblock {\em Statistical Science}, 1(3):297--310, 1986.

\bibitem{hastie2017extended}
Trevor Hastie, Robert Tibshirani, and Ryan~J Tibshirani.
\newblock Extended comparisons of best subset selection, forward stepwise
  selection, and the lasso.
\newblock Technical report, arXiv preprint 1707.08692, 2017.

\bibitem{herbrich2000}
R.~Herbrich, T.~Graepel, and K.~Obermayer.
\newblock {\em Large margin rank boundaries forordinal regression}.
\newblock MIT Press, 2000.

\bibitem{herbrich2001learning}
Ralf Herbrich.
\newblock {\em Learning Kernel Classifiers: Theory and Algorithms}.
\newblock MIT Press, 2001.

\bibitem{hornik1991approximation}
Kurt Hornik.
\newblock Approximation capabilities of multilayer feedforward networks.
\newblock {\em Neural Networks}, 4(2):251--257, 1991.

\bibitem{HYAFIL197615}
Laurent Hyafil and Ronald~L. Rivest.
\newblock Constructing optimal binary decision trees is {NP}-complete.
\newblock {\em Information Processing Letters}, 5(1):15--17, 1976.

\bibitem{icartetraining}
Rodrigo~Toro Icarte, Le{\'o}n Illanes, Margarita~P Castro, Andre~A Cire,
  Sheila~A McIlraith, and J~Christopher Beck.
\newblock Training binarized neural networks using {MIP} and {CP}.
\newblock In {\em Proceedings of the International Conference on Principles and
  Practice of Constraint Programming}, 2019.

\bibitem{izenman2008modern}
Alan~Julian Izenman.
\newblock {\em Modern Multivariate Statistical Techniques: Regression,
  Classification and Manifold Learning}, volume~10 of {\em Springer Texts in
  Statistics}.
\newblock Springer, 2008.

\bibitem{jaakkola2010learning}
Tommi Jaakkola, David Sontag, Amir Globerson, and Marina Meila.
\newblock Learning {B}ayesian network structure using {LP} relaxations.
\newblock In {\em Proceedings of the International Conference on Artificial
  Intelligence and Statistics}, pages 358--365, 2010.

\bibitem{jain1999data}
Anil~K. Jain, M.~N.~Narasimha Murty, and Patrick~J. Flynn.
\newblock Data clustering: a review.
\newblock {\em ACM Computing Surveys}, 31(3):264--323, 1999.

\bibitem{james2013introduction}
Gareth James, Daniela Witten, Trevor Hastie, and Robert Tibshirani.
\newblock {\em An Introduction to Statistical Learning}, volume 112.
\newblock Springer, 2013.

\bibitem{jan1994nonlinear}
Rong-Hong Jan and Maw-Sheng Chern.
\newblock Nonlinear integer bilevel programming.
\newblock {\em European Journal of Operational Research}, 72(3):574--587, 1994.

\bibitem{jolliffe2011principal}
Ian Jolliffe.
\newblock Principal component analysis.
\newblock In {\em International Encyclopedia of Statistical Science}, pages
  1094--1096. Springer, 2011.

\bibitem{karmitsa2017new}
Napsu Karmitsa, Adil~M. Bagirov, and Sona Taheri.
\newblock New diagonal bundle method for clustering problems in large data
  sets.
\newblock {\em European Journal of Operational Research}, 263(2):367--379,
  2017.

\bibitem{karpathy2014large}
Andrej Karpathy, George Toderici, Sanketh Shetty, Thomas Leung, Rahul
  Sukthankar, and Li~Fei-Fei.
\newblock Large-scale video classification with convolutional neural networks.
\newblock In {\em Proceedings of the IEEE conference on Computer Vision and
  Pattern Recognition}, pages 1725--1732, 2014.

\bibitem{katz2017ReLUplex}
Guy Katz, Clark Barrett, David~L. Dill, Kyle Julian, and Mykel~J. Kochenderfer.
\newblock Reluplex: An efficient {SMT} solver for verifying deep neural
  networks.
\newblock In {\em Proceedings of the International Conference on Computer Aided
  Verification}, pages 97--117, 2017.

\bibitem{kawano2015sparse}
Shuichi Kawano, Hironori Fujisawa, Toyoyuki Takada, and Toshihiko Shiroishi.
\newblock Sparse principal component regression with adaptive loading.
\newblock {\em Computational Statistics \& Data Analysis}, 89:192--203, 2015.

\bibitem{kelley1999iterative}
Carl~T Kelley.
\newblock {\em Iterative Methods for Optimization}.
\newblock Society for Industrial and Applied Mathematics, 1999.

\bibitem{KESHVARI2018585}
Abolfazl Keshvari.
\newblock Segmented concave least squares: A nonparametric piecewise linear
  regression.
\newblock {\em European Journal of Operational Research}, 266(2):585--594,
  2018.

\bibitem{khalil2016learning}
Elias~B. Khalil, Pierre~Le Bodic, Le~Song, George Nemhauser, and Bistra
  Dilkina.
\newblock Learning to branch in mixed integer programming.
\newblock In {\em Proceedings of the AAAI Conference on Artificial
  Intelligence}, pages 724--731, 2016.

\bibitem{khalil2017learning}
Elias~B. Khalil, Bistra Dilkina, George~L. Nemhauser, Shabbir Ahmed, and Yufen
  Shao.
\newblock Learning to run heuristics in tree search.
\newblock In {\em Proceedings of the International Joint Conference on
  Artificial Intelligence}, pages 659--666, 2017.

\bibitem{khalil2018combinatorial}
Elias~Boutros Khalil, Amrita Gupta, and Bistra Dilkina.
\newblock Combinatorial attacks on binarized neural networks.
\newblock Technical report, arXiv preprint 1810.03538, 2018.

\bibitem{Adam}
Diederik~P. Kingma and Jimmy Ba.
\newblock Adam: {A} method for stochastic optimization.
\newblock Technical report, arXiv preprint 1412.6980, 2014.

\bibitem{2017arXiv170207790H}
Diego Klabjan and Mark Harmon.
\newblock Activation ensembles for deep neural networks.
\newblock In {\em Proceeding of the IEEE International Conference on Big Data},
  pages 206--214, 2019.

\bibitem{klastorin1985p}
Ted~D Klastorin.
\newblock The p-median problem for cluster analysis: A comparative test using
  the mixture model approach.
\newblock {\em Management Science}, 31(1):84--95, 1985.

\bibitem{klatzer2015continuous}
Teresa Klatzer and Thomas Pock.
\newblock Continuous hyper-parameter learning for support vector machines.
\newblock In {\em Proceedings of the Computer Vision Winter Workshop}, pages
  39--47, 2015.

\bibitem{kramer2001}
Stefan Kramer, Gerhard Widmer, Bernhard Pfahringer, and Michael De~Groeve.
\newblock Prediction of ordinal classes using regression trees.
\newblock {\em Fundamenta Informaticae}, 47(1-2):1--13, 2001.

\bibitem{kraus2020deep}
Mathias Kraus, Stefan Feuerriegel, and Asil Oztekin.
\newblock Deep learning in business analytics and operations research: models,
  applications and managerial implications.
\newblock {\em European Journal of Operational Research}, 281(3):628--641,
  2020.

\bibitem{krizhevsky2012imagenet}
Alex Krizhevsky, Ilya Sutskever, and Geoffrey~E Hinton.
\newblock Imagenet classification with deep convolutional neural networks.
\newblock In F.~Pereira, C.~J.~C. Burges, L.~Bottou, and K.~Q. Weinberger,
  editors, {\em Advances in Neural Information Processing Systems}, pages
  1097--1105. Curran Associates, Inc., 2012.

\bibitem{kurakin2016adversarial}
Alexey Kurakin, Ian Goodfellow, and Samy Bengio.
\newblock Adversarial machine learning at scale.
\newblock Technical report, arXiv preprint 1611.01236, 2016.

\bibitem{kwatera1993clustering}
Renata~Krystyna Kwatera and Bruno Simeone.
\newblock Clustering heuristics for set covering.
\newblock {\em Annals of Operations Research}, 43(5):295--308, 1993.

\bibitem{lanckriet2002robust}
Gert R.~G. Lanckriet, Laurent~El Ghaoui, Chiranjib Bhattacharyya, and
  Michael~I. Jordan.
\newblock A robust minimax approach to classification.
\newblock {\em Journal of Machine Learning Research}, 3:555--582, 2002.

\bibitem{lecun_energy}
Yann LeCun, Sumit Chopra, Raia Hadsell, {Marc Aurelio} Ranzato, and {Fu Jie}
  Huang.
\newblock {\em A tutorial on energy-based learning}.
\newblock MIT Press, 2006.

\bibitem{lecun1989generalization}
Yann LeCun et~al.
\newblock Generalization and network design strategies.
\newblock {\em Connectionism in Perspective}, 19:143--155, 1989.

\bibitem{leof2018automver}
Francesco Leofante, Nina Narodytska, Luca Pulina, and Armando Tacchella.
\newblock Automated verification of neural networks: Advances, challenges and
  perspectives.
\newblock Technical report, arXiv preprint 1805.09938, 2018.

\bibitem{Lewis2014}
Mark Lewis, Haibo Wang, and Gary Kochenberger.
\newblock Exact solutions to the capacitated clustering problem: A comparison
  of two models.
\newblock {\em Annals of Data Science}, 1(1):15--23, 2014.

\bibitem{liang2017fisher}
Tengyuan Liang, Tomaso Poggio, Alexander Rakhlin, and James Stokes.
\newblock Fisher-{R}ao metric, geometry, and complexity of neural networks.
\newblock In {\em Proceeding of the International Conference on Artificial
  Intelligence and Statistics}, pages 888--896, 2019.

\bibitem{lin2013network}
Min Lin, Qiang Chen, and Shuicheng Yan.
\newblock Network in network.
\newblock Technical report, arXiv preprint 1312.4400, 2013.

\bibitem{liu2016teaching}
Ji~Liu and Xiaojin Zhu.
\newblock The teaching dimension of linear learners.
\newblock {\em The Journal of Machine Learning Research}, 17(1):5631--5655,
  2016.

\bibitem{LodZ17}
Andrea Lodi and Giulia Zarpellon.
\newblock On learning and branching: a survey.
\newblock {\em TOP}, 25(2):207--236, 2017.

\bibitem{lombardi2017empirical}
Michele Lombardi, Michela Milano, and Andrea Bartolini.
\newblock Empirical decision model learning.
\newblock {\em Artificial Intelligence}, 244:343--367, 2017.

\bibitem{lowd2005adversarial}
Daniel Lowd and Christopher Meek.
\newblock Adversarial learning.
\newblock In {\em Proceedings of the International Conference on Knowledge
  Discovery in Data Mining}, pages 641--647, 2005.

\bibitem{macqueen1967some}
James MacQueen.
\newblock Some methods for classification and analysis of multivariate
  observations.
\newblock In {\em Proceedings of the Berkeley Symposium on Mathematical
  Statistics and Probability}, pages 281--297, 1967.

\bibitem{madry2017towards}
Aleksander Madry, Aleksandar Makelov, Ludwig Schmidt, Dimitris Tsipras, and
  Adrian Vladu.
\newblock Towards deep learning models resistant to adversarial attacks.
\newblock Technical report, arXiv preprint 1706.06083, 2017.

\bibitem{MAI2018594}
Feng Mai, Michael~J. Fry, and Jeffrey~W. Ohlmann.
\newblock Model-based capacitated clustering with posterior regularization.
\newblock {\em European Journal of Operational Research}, 271(2):594--605,
  2018.

\bibitem{maldonado2014feature}
Sebasti{\'a}n Maldonado, Juan P{\'e}rez, Richard Weber, and Martine Labb{\'e}.
\newblock Feature selection for support vector machines via mixed integer
  linear programming.
\newblock {\em Information Sciences}, 279:163--175, 2014.

\bibitem{mei2015using}
Shike Mei and Xiaojin Zhu.
\newblock Using machine teaching to identify optimal training-set attacks on
  machine learners.
\newblock In {\em Proceedings of the AAAI Conference on Artificial
  Intelligence}, pages 2871--2877, 2015.

\bibitem{mielke1997permutation}
Paul~W Mielke and Kenneth~J Berry.
\newblock Permutation-based multivariate regression analysis: The case for
  least sum of absolute deviations regression.
\newblock {\em Annals of Operations Research}, 74:259, 1997.

\bibitem{miller2002subset}
Alan Miller.
\newblock {\em Subset Selection in Regression}.
\newblock Chapman and Hall/CRC, 2002.

\bibitem{mivsic2017optimization}
Velibor~V Mi{\v{s}}i{\'c}.
\newblock Optimization of tree ensembles.
\newblock {\em Operations Research (In press)}, 2020.

\bibitem{miyashiro2015mixed}
Ryuhei Miyashiro and Yuichi Takano.
\newblock Mixed integer second-order cone programming formulations for variable
  selection in linear regression.
\newblock {\em European Journal of Operational Research}, 247(3):721--731,
  2015.

\bibitem{montufar2017notes}
Guido Mont{\'u}far.
\newblock Notes on the number of linear regions of deep neural networks.
\newblock Technical report, Max Planck Institute for Mathematics in the
  Sciences, 2017.

\bibitem{moore2011model}
Gregory Moore, Charles Bergeron, and Kristin~P Bennett.
\newblock Model selection for primal {SVM}.
\newblock {\em Machine Learning}, 85(1-2):175--208, 2011.

\bibitem{MORTENSON2015583}
Michael~J. Mortenson, Neil~F. Doherty, and Stewart Robinson.
\newblock Operational research from taylorism to terabytes: A research agenda
  for the analytics age.
\newblock {\em European Journal of Operational Research}, 241(3):583--595,
  2015.

\bibitem{mulvey1979cluster}
John~M Mulvey and Harlan~P Crowder.
\newblock Cluster analysis: An application of {L}agrangian relaxation.
\newblock {\em Management Science}, 25(4):329--340, 1979.

\bibitem{natarajan1995sparse}
Balas~Kausik Natarajan.
\newblock Sparse approximate solutions to linear systems.
\newblock {\em SIAM Journal on Computing}, 24(2):227--234, 1995.

\bibitem{negreiros2006capacitated}
Marcos Negreiros and Augusto Palhano.
\newblock The capacitated centred clustering problem.
\newblock {\em Computers \& Operations Research}, 33(6):1639--1663, 2006.

\bibitem{nie2015learning}
Siqi Nie, Cassio~P De~Campos, and Qiang Ji.
\newblock Learning bounded tree-width {B}ayesian networks via sampling.
\newblock In {\em Proceedings of the European Conference on Symbolic and
  Quantitative Approaches to Reasoning and Uncertainty}, pages 387--396.
  Springer, 2015.

\bibitem{nie2014advances}
Siqi Nie, Denis~D Mau{\'a}, Cassio~P De~Campos, and Qiang Ji.
\newblock Advances in learning {B}ayesian networks of bounded treewidth.
\newblock In {\em Advances in Neural Information Processing Systems}, pages
  2285--2293, 2014.

\bibitem{OLAFSSON20081429}
Sigurdur Olafsson, Xiaonan Li, and Shuning Wu.
\newblock Operations research and data mining.
\newblock {\em European Journal of Operational Research}, 187(3):1429 -- 1448,
  2008.

\bibitem{parviainen2014learning}
Pekka Parviainen, Hossein~Shahrabi Farahani, and Jens Lagergren.
\newblock Learning bounded tree-width {B}ayesian networks using integer linear
  programming.
\newblock In {\em Artificial Intelligence and Statistics}, pages 751--759,
  2014.

\bibitem{patil2014optimal}
Kaustubh~R Patil, Jerry Zhu, \L~ukasz Kope\'{c}, and Bradley~C Love.
\newblock Optimal teaching for limited-capacity human learners.
\newblock In Z.~Ghahramani, M.~Welling, C.~Cortes, N.~D. Lawrence, and K.~Q.
  Weinberger, editors, {\em Advances in Neural Information Processing Systems
  27}, pages 2465--2473. Curran Associates, Inc., 2014.

\bibitem{1674938}
Harold~J. Payne and William~S. Meisel.
\newblock An algorithm for constructing optimal binary decision trees.
\newblock {\em IEEE Transactions on Computers}, 26(9):905--916, 1977.

\bibitem{scikit-learn}
Fabian Pedregosa, Ga{\"e}l Varoquaux, Alexandre Gramfort, Vincent Michel,
  Bertrand Thirion, Olivier Grisel, Mathieu Blondel, Peter Prettenhofer, Ron
  Weiss, Vincent Dubourg, Alexandre~Passos Jake Vanderpla~and, David
  Cournapeau, Matthieu Brucher, Matthieu Perrot, and Édouard Duchesnay.
\newblock Scikit-learn: Machine learning in {P}ython.
\newblock {\em Journal of Machine Learning Research}, 12:2825--2830, 2011.

\bibitem{PIRAMUTHU2004483}
Selwyn Piramuthu.
\newblock Evaluating feature selection methods for learning in data mining
  applications.
\newblock {\em European Journal of Operational Research}, 156(2):483--494,
  2004.

\bibitem{poggio2017theory}
Tomaso Poggio, Kenji Kawaguchi, Qianli Liao, Brando Miranda, Lorenzo Rosasco,
  Xavier Boix, Jack Hidary, and Hrushikesh Mhaskar.
\newblock Theory of deep learning {III}: explaining the non-overfitting puzzle.
\newblock Technical report, arXiv preprint 1801.00173, 2017.

\bibitem{reris2015principal}
Robert Reris and Jean~Paul Brooks.
\newblock Principal component analysis and optimization: a tutorial.
\newblock Technical report, Virginia Commonwealth University, 2015.

\bibitem{robbins1985stochastic}
Herbert Robbins and Sutton Monro.
\newblock A stochastic approximation method.
\newblock In {\em Herbert Robbins Selected Papers}, pages 102--109. Springer,
  1985.

\bibitem{rovatti2014optimistic}
Riccardo Rovatti, Claudia D'Ambrosio, Andrea Lodi, and Silvano Martello.
\newblock Optimistic {MILP} modeling of non-linear optimization problems.
\newblock {\em European Journal of Operational Research}, 239(1):32--45, 2014.

\bibitem{sauglam2006mixed}
Burcu Sa{\u{g}}lam, F.~Sibel Salman, Serpil Say{\i}n, and Metin T{\"u}rkay.
\newblock A mixed-integer programming approach to the clustering problem with
  an application in customer segmentation.
\newblock {\em European Journal of Operational Research}, 173(3):866--879,
  2006.

\bibitem{santi2016model}
{\'E}verton Santi, Daniel Aloise, and Simon~J. Blanchard.
\newblock A model for clustering data from heterogeneous dissimilarities.
\newblock {\em European Journal of Operational Research}, 253(3):659--672,
  2016.

\bibitem{scanagatta2016learning}
Mauro Scanagatta, Giorgio Corani, Cassio~P de~Campos, and Marco Zaffalon.
\newblock Learning treewidth-bounded bayesian networks with thousands of
  variables.
\newblock In D.~D. Lee, M.~Sugiyama, U.~V. Luxburg, I.~Guyon, and R.~Garnett,
  editors, {\em Advances in Neural Information Processing Systems}, pages
  1462--1470. Curran Associates, Inc., 2016.

\bibitem{SCHEUERER2006533}
Stephan Scheuerer and Rolf Wendolsky.
\newblock A scatter search heuristic for the capacitated clustering problem.
\newblock {\em European Journal of Operational Research}, 169(2):533--547,
  2006.

\bibitem{schobel1998locating}
Anita Sch{\"o}bel.
\newblock Locating least-distant lines in the plane.
\newblock {\em European Journal of Operational Research}, 106(1):152--159,
  1998.

\bibitem{serra2017bounding}
Thiago Serra, Christian Tjandraatmadja, and Srikumar Ramalingam.
\newblock Bounding and counting linear regions of deep neural networks.
\newblock In {\em Proceeding of the International Conference on Machine
  Learning}, pages 4558--4566, 2018.

\bibitem{shaham2018understanding}
Uri Shaham, Yutaro Yamada, and Sahand Negahban.
\newblock Understanding adversarial training: Increasing local stability of
  supervised models through robust optimization.
\newblock {\em Neurocomputing}, 307:195--204, 2018.

\bibitem{shashua2003ranking}
Amnon Shashua and Anat Levin.
\newblock Ranking with large margin principle: Two approaches.
\newblock In S.~Becker, S.~Thrun, and K.~Obermayer, editors, {\em Advances in
  Neural Information Processing Systems}, pages 961--968. MIT Press, 2003.

\bibitem{Shinohara1991}
Ayumi Shinohara and Satoru Miyano.
\newblock Teachability in computational learning.
\newblock {\em New Generation Computing}, 8(4):337--347, Feb 1991.

\bibitem{smola2004tutorial}
Alex~J Smola and Bernhard Sch{\"o}lkopf.
\newblock A tutorial on support vector regression.
\newblock {\em Statistics and Computing}, 14(3):199--222, 2004.

\bibitem{solomonoff1957inductive}
Raymond~J. Solomonoff.
\newblock An inductive inference machine.
\newblock In {\em IRE Convention Record, Section on Information Theory},
  volume~2, pages 56--62, 1957.

\bibitem{song2019review}
Heda Song, Isaac Triguero, and Ender {\"O}zcan.
\newblock A review on the self and dual interactions between machine learning
  and optimisation.
\newblock {\em Progress in Artificial Intelligence}, 8(2):143--165, 2019.

\bibitem{SteinhardtKL17}
Jacob Steinhardt, Pang Wei~W Koh, and Percy~S Liang.
\newblock Certified defenses for data poisoning attacks.
\newblock In I.~Guyon, U.~V. Luxburg, S.~Bengio, H.~Wallach, R.~Fergus,
  S.~Vishwanathan, and R.~Garnett, editors, {\em Advances in Neural Information
  Processing Systems}, pages 3517--3529. Curran Associates, Inc., 2017.

\bibitem{szegedy2013intriguing}
Christian Szegedy, Wojciech Zaremba, Ilya Sutskever, Joan Bruna, Dumitru Erhan,
  Ian Goodfellow, and Rob Fergus.
\newblock Intriguing properties of neural networks.
\newblock Technical report, arXiv preprint 1312.6199, 2013.

\bibitem{tamura2017best}
Ryuta Tamura, Ken Kobayashi, Yuichi Takano, Ryuhei Miyashiro, Kazuhide Nakata,
  and Tomomi Matsui.
\newblock Best subset selection for eliminating multicollinearity.
\newblock {\em Journal of the Operations Research Society of Japan},
  60(3):321--336, 2017.

\bibitem{tamura2019mixed}
Ryuta Tamura, Ken Kobayashi, Yuichi Takano, Ryuhei Miyashiro, Kazuhide Nakata,
  and Tomomi Matsui.
\newblock Mixed integer quadratic optimization formulations for eliminating
  multicollinearity based on variance inflation factor.
\newblock {\em Journal of Global Optimization}, 73(2):431--446, 2019.

\bibitem{taylan2007new}
Pakize Taylan, G-W Weber, and Amir Beck.
\newblock New approaches to regression by generalized additive models and
  continuous optimization for modern applications in finance, science and
  technology.
\newblock {\em Optimization}, 56(5-6):675--698, 2007.

\bibitem{TIWARI2018319}
Sunil Tiwari, H.M. Wee, and Yosef Daryanto.
\newblock Big data analytics in supply chain management between 2010 and 2016:
  Insights to industries.
\newblock {\em Computers \& Industrial Engineering}, 115:319--330, 2018.

\bibitem{tjeng2017verifying}
Vincent Tjeng and Russ Tedrake.
\newblock Evaluating robustness of neural networks with mixed integer
  programming.
\newblock Technical report, arXiv preprint 1711.07356, 2017.

\bibitem{TORIELLO201286}
Alejandro Toriello and Juan~Pablo Vielma.
\newblock Fitting piecewise linear continuous functions.
\newblock {\em European Journal of Operational Research}, 219(1):86--95, 2012.

\bibitem{tramer2017ensemble}
Florian Tram{\`e}r, Alexey Kurakin, Nicolas Papernot, Ian Goodfellow, Dan
  Boneh, and Patrick McDaniel.
\newblock Ensemble adversarial training: Attacks and defenses.
\newblock Technical report, arXiv preprint 1705.07204, 2017.

\bibitem{vapnik1998statistical}
Vladimir Vapnik.
\newblock {\em Statistical Learning Theory}, volume~3.
\newblock Wiley, New York, 1998.

\bibitem{vapnik2013nature}
Vladimir Vapnik.
\newblock {\em The Nature of Statistical Learning Theory}.
\newblock Springer Science \& Business Media, 2013.

\bibitem{verwer2017learning}
Sicco Verwer and Yingqian Zhang.
\newblock Learning decision trees with flexible constraints and objectives
  using integer optimization.
\newblock In {\em Proceedings of the International Conference on AI and OR
  Techniques in Constraint Programming for Combinatorial Optimization
  Problems}, pages 94--103, 2017.

\bibitem{verwer2017auction}
Sicco Verwer, Yingqian Zhang, and Qing~Chuan Ye.
\newblock Auction optimization using regression trees and linear models as
  integer programs.
\newblock {\em Artificial Intelligence}, 244:368--395, 2017.

\bibitem{vielma2010mixed}
Juan~Pablo Vielma, Shabbir Ahmed, and George Nemhauser.
\newblock Mixed-integer models for nonseparable piecewise-linear optimization:
  Unifying framework and extensions.
\newblock {\em Operations Research}, 58(2):303--315, 2010.

\bibitem{VACLAVIK2018}
Roman Václavík, Antonín Novák, Přemysl Šůcha, and Zdeněk Hanzálek.
\newblock Accelerating the branch-and-price algorithm using machine learning.
\newblock {\em European Journal of Operational Research}, 271(3):1055--1069,
  2018.

\bibitem{WANG201698}
Gang Wang, Angappa Gunasekaran, Eric~W.T. Ngai, and Thanos Papadopoulos.
\newblock Big data analytics in logistics and supply chain management: Certain
  investigations for research and applications.
\newblock {\em International Journal of Production Economics}, 176:98--110,
  2016.

\bibitem{wang2010multi}
Hua Wang, Chris Ding, and Heng Huang.
\newblock Multi-label linear discriminant analysis.
\newblock In {\em Proceedings of the European Conference on Computer Vision},
  pages 126--139, 2010.

\bibitem{wang2006doubly}
Li~Wang, Ji~Zhu, and Hui Zou.
\newblock The doubly regularized support vector machine.
\newblock {\em Statistica Sinica}, 16(2):589, 2006.

\bibitem{WangPoisoning18}
Yizhen Wang and Kamalika Chaudhuri.
\newblock Data poisoning attacks against online learning.
\newblock Technical report, arXiv preprint 1808.08994, 2018.

\bibitem{wang2018enhancing}
Yuan Wang, Dongxiang Zhang, Ying Liu, Bo~Dai, and Loo~Hay Lee.
\newblock Enhancing transportation systems via deep learning: A survey.
\newblock {\em Transportation Research Part C: Emerging Technologies},
  99:144--163, 2019.

\bibitem{wistuba2019survey}
Martin Wistuba, Ambrish Rawat, and Tejaswini Pedapati.
\newblock A survey on neural architecture search.
\newblock Technical report, arXiv preprint arXiv:1905.01392, 2019.

\bibitem{wright2017optimization}
Stephen~J Wright.
\newblock Optimization algorithms for data analysis.
\newblock In Michael Mahoney, John Duchi, and Anna Gilbert, editors, {\em The
  Mathematics of Data}, pages 49--98. American Mathematical Society, 2018.

\bibitem{yuan2013learning}
Changhe Yuan and Brandon Malone.
\newblock Learning optimal {B}ayesian networks: A shortest path perspective.
\newblock {\em Journal of Artificial Intelligence Research}, 48:23--65, 2013.

\bibitem{zhu2013machine}
Jerry Zhu.
\newblock Machine teaching for bayesian learners in the exponential family.
\newblock In C.~J.~C. Burges, L.~Bottou, M.~Welling, Z.~Ghahramani, and K.~Q.
  Weinberger, editors, {\em Advances in Neural Information Processing Systems},
  pages 1905--1913. Curran Associates, Inc., 2013.

\bibitem{zhu20041}
Ji~Zhu, Saharon Rosset, Robert Tibshirani, and Trevor~J. Hastie.
\newblock 1-{N}orm support vector machines.
\newblock In S.~Thrun, L.~K. Saul, and B.~Sch\"{o}lkopf, editors, {\em Advances
  in Neural Information Processing Systems}, pages 49--56. MIT Press, 2004.

\bibitem{zhu2015machine}
Xiaojin Zhu.
\newblock Machine teaching: An inverse problem to machine learning and an
  approach toward optimal education.
\newblock In {\em Proceedings of the AAAI Conference on Artificial
  Intelligence}, pages 4083--4087, 2015.

\bibitem{zhu2018overview}
Xiaojin Zhu, Adish Singla, Sandra Zilles, and Anna~N Rafferty.
\newblock An overview of machine teaching.
\newblock Technical report, arXiv preprint 1801.05927, 2018.

\bibitem{zinkevich2003online}
Martin Zinkevich.
\newblock Online convex programming and generalized infinitesimal gradient
  ascent.
\newblock In {\em Proceedings of the International Conference on Machine
  Learning}, pages 928--936, 2003.

\bibitem{zou2005regularization}
Hui Zou and Trevor Hastie.
\newblock Regularization and variable selection via the elastic net.
\newblock {\em Journal of the Royal Statistical Society: Series B (Statistical
  Methodology)}, 67(2):301--320, 2005.

\end{thebibliography}
\newpage

\end{document}